\documentclass[12pt]{extarticle}
\usepackage{amsmath, amsthm, amssymb, hyperref, color}
\usepackage{graphicx}
\usepackage{caption}
\usepackage{subcaption}
\usepackage{mathtools}
\usepackage{enumerate}
\usepackage{stackengine}
\usepackage{verbatim}
\usepackage{upgreek, adjustbox}
\usepackage{chemfig, chemnum}
\usepackage{tikz}
\usetikzlibrary{patterns}
\usepackage[lined,commentsnumbered,ruled]{algorithm2e}
\usepackage{caption}
\usepackage{subcaption}
\usepackage{float}
\oddsidemargin \evensidemargin
\usetikzlibrary{decorations.pathreplacing}
\newtheorem{theorem}{Theorem}[section]

\theoremstyle{definition}
\newtheorem{definition}[theorem]{Definition}

\newtheorem{remark}[theorem]{Remark}
\newtheorem{conjecture}[theorem]{Conjecture}

\newtheorem{example}[theorem]{Example}

\newcommand{\CC}{\mathbb{C} }
\newcommand{\ZZ}{\mathbb{Z}}

\usepackage{color}
\newcommand{\arrowIn}{\tikz \draw[-stealth] (-1pt,0) -- (1pt,0);}

\title{\textbf{Crossing the transcendental divide: from translation surfaces to algebraic curves}}

\author{{T\"urk\"u \"Ozl\"um \c{C}elik$^1$, Samantha Fairchild$^{2, 3}$, Yelena Mandelshtam$^4$}}
\date{\footnotesize{$^1$Bo\u{g}azi\c{c}i University, $^2$University of Osnabr\"uck, $^3$MPI-MiS Leipzig, $^4$UC Berkeley}}

\begin{document}

\maketitle
%\tableofcontents

\noindent \begin{abstract}
 We study constructing an algebraic curve from a Riemann surface given via a translation surface, which is a collection of finitely many polygons in the plane with sides identified by translation. We use the theory of discrete Riemann surfaces to give an algorithm for approximating the Jacobian variety of a translation surface whose polygon can be decomposed into squares. We first implement the algorithm in the case of $L$ shaped polygons where the algebraic curve is already known. The algorithm is also implemented in any genus for specific examples of Jenkins--Strebel representatives, a dense family of translation surfaces that, until now, lived squarely on the analytic side of the transcendental divide between Riemann surfaces and algebraic curves. Using Riemann theta functions, we give numerical experiments and resulting conjectures up to genus $5$.
\end{abstract}

\section{Introduction}
We present an algorithm with numerical experiments as a step in bridging the \emph{transcendental divide} between Riemann surfaces and algebraic curves. The classical equivalence of Riemann surfaces and algebraic curves leaves a divide in the sense that it is non-trivial to determine the exact curve associated to a Riemann surface. This becomes transcendental as connecting a Riemann surface to an algebraic curve utilizes Riemann theta functions. In our case we use discrete harmonic functions to approximate the Jacobian variety of the algebraic curve of a translation surface and the theta functions to understand these approximations.%~\cite{Mum} 

%%%%%%%%\textbf{Goals:} 

 A translation surface is obtained by identifying edges of finitely many polygons in the plane with  complex translations. This identification of polygons builds a Riemann surface $X$, or equivalently a complex algebraic curve. The surface $X$ naturally comes equipped with a nonzero holomorphic $1$-form $\omega$, given locally by $dz$, with $2g-2$ zeroes located at the vertices of the polygons, where $g$ is the genus of $X$. One method to obtain an equation for the underlying algebraic curve of a translation surface is to first find a basis of the space of holomorphic one-forms on $X$. This space identifies a canonical model of the curve in some projective space. For instance, if the translation surface has extra automorphisms the basis of holomorphic one-forms can be determined exactly \cite{DM21, Rod, Sil}. 
Instead of requiring extra symmetries, we here aim to construct an algebraic curve explicitly from a translation surface by approximating the Riemann matrix.

Recall that a Riemann matrix $\tau$ associated to an algebraic curve $\mathcal{C}$ is defined by integrating a canonical basis of holomorphic differentials over the cycles forming a homology basis for the curve. The matrix $\tau$ data defines the Jacobian variety of $\mathcal{C}$, namely the quotient $\CC^g/(\ZZ^g+\tau\ZZ^g)$.
The Torelli theorem states that the curve $\mathcal{C}$ is determined by its Jacobian variety. In practice, one can construct Riemann matrices given an algebraic curve \cite{Decon}, and, less trivially, given $\tau$ reconstruct the algebraic curve at least in low genus \cite{AgoCelEke}. We will approximate the Riemann matrix from a translation surface and then utilize tools from numerical algebraic geometry to approximate the underlying curve.

 %%%%%%\textbf{Tools.} 
 
 The key instruments of our approach are discrete Riemann surfaces. %(\S~\ref{sec:DiscreteApproximations})
They can be thought as discrete counterparts of the Riemann surfaces in discrete complex analysis~\cite{Smir}. All notions such as Riemann matrices and holomorphic differential forms have corresponding discretizations. An important result in the literature of discrete Riemann surfaces, which is fundamental for the present study, uses discrete energies to prove convergence of the discrete Riemann surfaces to the underlying Riemann surface. In particular the discrete Riemann matrices converge to a Riemann matrix~\cite{BobBuc, BobSko} of the Riemann surface. Our experiments employ compact Riemann surfaces given by quadrangulations based on theory developed in \cite{BobGun}. In the case of \cite{BobGun}, convergence of the discrete Riemann matrices to a Riemann surface is not known in full generality, but our results give evidence for a potential result on convergence of the period matrices of \cite{BobGun}.
 
 %The convergence rate is given by looking at particular triangulations as the longest edge length $h$ goes to zero. For translation surfaces, the convergence rate \cite[Theorem 2.5]{BobSko}, is proportional to $h$ when $g=1$. When $g\geq 2$, the rate is $h|\log(h)|$ when the $2g-2$ zeroes are all distinct, and at worst $h^\frac{2}{2g-1}$ which corresponds to the case when there is a single zero of order $2g-2$. In \cite[Theorem 2.4]{BobBuc}, when working with ramified coverings of the sphere and carefully choosing an optimal triangulation, the convergence rate is always proportional to $h$. We consider square quadrangulations instead of triangulations. After discussions with Felix G{\"u}nther, we expect similar convergence rates to \cite{BobSko}. The key idea is that by dividing the squares along the diagonal, the resulting triangulation is sufficiently regular to fit into the framework of \cite{BobSko}. 

%%%%%%%%%%%%%\textbf{Algorithms.} 

A fairly technical algorithm to compute discrete Riemann matrices when the discrete complex structure is given by triangulations was presented in \cite{Bob11}. In contrast to the previous work, we present Algorithm~\ref{alg:general} to be accessible to a large variety of mathematical audiences. Algorithm~\ref{alg:general} inputs a translation surface as a polygon, a level of approximation, and outputs the associated discrete Riemann matrix. %The technical details that a person implementing the algorithm must know are as follows. 
The given polygon must be able to be divided into squares. Moreover one must place an initial bipartite graph on the square tiled polygons, which means that the vertices are either black and white, and no two vertices of the same color are connected by an edge. The bipartite graph must be chosen so that all identified vertices of the translation surface are the same color. The rest is to find the correct basis of homology which respects the pairings given in the translation surface.

We present two concrete implementations of Algorithm~\ref{alg:general}. In Algorithm~\ref{alg:L} we consider a family of symmetric $L$ shaped polygons which are all genus $2$ Riemann surfaces. The $L$ shapes are highly symmetric, so they serve as a good test case where the exact underlying algebraic curve is known~\cite{Sil}. %\S~\ref{sec:Translation}. 
Another benefit of the $L$ shapes is that they allow us to experiment with convergence to polygons which cannot be square tiled e.g., shapes with an irrational side length, by closer approximations of square tiled polygons. 

We also present Algorithm~\ref{sec:AlgForJS}, which gives a natural family of square tiled translation surfaces, called Jenkins--Strebel differentials, for any genus $g\geq 2$. When working with any genus, the step of Algorithm~\ref{alg:general} where we must choose a homology basis respecting the identifications requires care. We explain our difficulties, and explain how we overcame these difficulties. This algorithm leads to two interesting experiments. First we approximate the Riemann matrix in a case where the underlying curve is not yet known. In genus 2, we can numerically compute the hyperelliptic curve, leading to some conjectures on the structure of the underlying curves. Further, we can do experiments in genus $g= 3, 4,$ and $5$, to understand how the discrete matrices approach along the Schottky locus of $g\times g$ matrices which are associated to an algebraic curve. 

In Section~\ref{sec:background}, we give history, definitions, and examples for algebraic curves, translation surfaces, and discrete Riemann surfaces. In Section~\ref{sec:algorithms}, we present the algorithms and numerical experiments introduced above. Namely in Algorithm~\ref{alg:general} we give the algorithm for any translation surface. Given a translation surface via its defining polygons, we approximate its Riemann matrix through a family of discrete Riemann surfaces, which are expressed in terms of subdivisions of the polygons. We then run experiments in two specific cases. In Section~\ref{sec:Lsection} we describe the algorithm in details and discuss convergence to the underlying Riemann surface for the case of the $L$. In Section~\ref{sec:JS} we again give specifics of the algorithm for the family of Jenkins-Strebel differentials, and we use theta functions to approximate an equation of the algebraic curve from the estimated matrix. Finally in Appendix~\ref{appendix} we give tables of values for approximating discrete Riemann matrices. The code can be found at our MathRepo page~\cite{CelFaiMan}.
\newline
\newline\noindent{\bf Acknowledgements.} We wish to thank Felix G{\"u}nther for sharing his knowlege of discrete Riemann surfaces, Bernd Sturmfels for connecting us on this project, and the anonymous referee for detailed and constructive feedback. We also want to thank Nils Bruin, Christophe Ritzenthaler, and Andr\'{e} Uschmajew  for useful conversations. T.\"O.{\c{C}} was supported by the Turkish Scientific and Technological Research Council (T\"{U}B{\.I}TAK) - Project number 121C039. S.F. was supported by the Deutsche Forschungsgemeinschaft
(DFG) – Projektnummer 445466444. Y.M. was partially supported by NSF grant DGE 2146752. The authors report there are no competing interests to declare.

\section{Background}\label{sec:background}

We aim here to recall some background which we need to reconstruct algebraic curves from their translation surfaces. This includes some preliminaries for Riemann surfaces and theta functions, translation surfaces, and discrete Riemann surfaces. Each topic is a well studied and interesting subject in its own right, so we provide references for their underlying theories.

\subsection{Riemann Surfaces, analytic and algebraic}\label{sec:RS}
%Theta functions introduction here

Riemann surfaces are one-dimensional complex manifolds, among which the compact ones are complex smooth algebraic curves. Among central objects underlying the connection between the analytic side and the algebraic side are the so-called \emph{theta functions}. 

Let $\mathcal{C}$ be a complex smooth algebraic curve. Let $\omega_1,\dots , \omega_g$ be a basis of ${\rm H}^0(\mathcal{C},\Omega^1_\mathcal{C})$, i.e. the space of holomorphic differential one-forms. Let $\alpha_1,\dots , \alpha_g, \beta_1, \dots, \beta_g$ be a symplectic basis of ${\rm H}_1(\mathcal{C},\mathbb{Z})$. The $g\times 2g$ matrix is called the \emph{period matrix}: 
\begin{equation}\label{eq:periodMatrix}
(\tau_1|\tau_2):=\left(\left(\int_{\alpha_i} \omega_j\right)\Big|\left(\int_{\beta_i} \omega_j\right)\right)
\end{equation}
%\end{Large}
and $\tau:=\tau_1^{-1}\tau_2$ is called a \emph{Riemann matrix} of the algebraic curve $\mathcal{C}$. The Riemann matrix lies in the set of $g \times g$ symmetric matrices with positive definite imaginary part with complex entries, the so-called Siegel upper half space $\mathbb{H}_g$. The theta function with characteristic $\varepsilon,\delta \in \{ 0,1\}^g$ is a complex-valued function defined on $\CC^g\times \mathbb{H}_g$: 
 \begin{equation}
\label{eq:thetaFunctionChar}
\theta\begin{bmatrix} \varepsilon \\ \delta \end{bmatrix}({\bf z}\, |\, \tau)\,\,\, = \,\,\,
\sum_{{\bf n} \in \mathbb{Z}^g} {\rm exp} \left( \pi \mathbf{i} \left({\bf n}+\frac{\varepsilon}{2}\right)^T \tau \left({\bf n}+\frac{\varepsilon}{2}\right) + \left({\bf n} + \frac{\varepsilon}{2}\right)^T  \left({\bf z} + \frac{\delta}{2}\right) \right).
\end{equation}
Here $\mathbf{i} = \sqrt{-1}$ to allow the use of the indexing variable $i$. When $\varepsilon=\delta=0$, this is nothing but the \emph{Riemann theta function}, and differs by an exponential factor from the latter. The characteristic is called even, odd if $\varepsilon\cdot \delta\equiv 0,1 \pmod 2$ respectively. So there are $2^{g-1}(2^g+1)$ odd and $2^{g-1}(2^g+1)$ even characteristics. For fixed $\tau$, the values $\theta\begin{bmatrix} \varepsilon \\ \delta \end{bmatrix}({\bf 0}\, |\, \tau)$ at ${\bf z}= {\bf 0}$ are known as \emph{theta constants}. We will also use the term theta constant for the evaluation of the derivatives of the theta function at ${\bf 0}$, which we denote as follows: 
\begin{equation}\label{thetaConstant}
  \theta^{\varepsilon,\delta}_{ij\dots}:=\frac{\partial}{\partial z_i}\frac{\partial}{\partial z_j}\dots\theta\begin{bmatrix} \varepsilon \\ \delta \end{bmatrix}({\bf z}\, |\, \tau)_{\big|_{{\bf z} ={\bf 0}}}.
\end{equation}

 The theta functions play a central role in the literature of the Schottky problem and the Torelli theorem~\cite{Gru12}. When $g=4$, the so-called Schotkky-Igusa modular form defines an analytic hypersurface~\cite[Theorem 1]{Igusa} in terms of theta functions, which describes Riemann matrices in the Siegel upper half space that are of algebraic curves, the so-called \emph{Schottky locus}. For higher genus, there are analytical equations in terms of theta functions defining a locus containing the Schottky locus. In the context of the Torelli theorem, we suggest the reader to see~\cite[Theorem 8.1]{Gua2002} for hyperelliptic curves and see~\cite[Theorem 1.1]{Gua2011} for genus 3 non-hyperelliptic curves. These rely on the classical formulae going back to Riemann, namely the Thomae formula and the Weber formula~\cite{Riemann, Thomae, Weber} with their generalizations to any genus~\cite{Cel2019}, which relate the extrinsic and intrinsic sides of geometry of the underlying curve. For instance, the vanishing theta constants i.e., the theta constants that are zero, are well-understood when the curve is hyperelliptic~\cite{Mum}. Actually, the theta constants express certain divisors of the curve $\mathcal{C}$, e.g. theta characteristic divisors (semi-canonical divisors), which recover the curve itself. More precisely, the choice of the basis $\omega_1,\ldots, \omega_g$ of the holomorphic differentials on $\mathcal{C}$ gives a map, namely the \emph{canonical map} of $\mathcal{C}$:
\begin{align*}
    \Phi :\,\, & \mathcal{C} \rightarrow \mathbb{P}^{g-1}\\
    & P \mapsto (\omega_1(P),\dots , \omega_g(P)).
\end{align*}
An important note is that the canonical image $\Phi(\mathcal{C})$ is defined over a field over which the differential forms are defined. The following statement formulates the theta characteristic divisors in terms of the canonical model: 

\begin{theorem}[Theorem 2.2, \cite{Gua2002}]
     Let $\tau_1$ be the first $g \times g$ part of the period matrix~\eqref{eq:periodMatrix}. Let $D$ be an effective theta characteristic divisor of degree $g-1$ with $\dim{\rm H}^0(\mathcal{C},D)=1$. The corresponding equations of the hyperplanes spanned by $\Phi(D)$ are given by:
    
    \begin{equation}\label{eq:thetaChar}
    \begin{pmatrix}
    \theta^{\varepsilon,\delta}_1 & ,\dots , & \theta^{\varepsilon,\delta}_g
    \end{pmatrix} \cdot \tau_1^{-1} \cdot \begin{pmatrix}
    X_1 \\ \vdots \\ X_g
    \end{pmatrix}
   = 0,
    \end{equation} where the characteristic $\begin{bmatrix} \varepsilon \\ \delta  \end{bmatrix}$
    ranges over the odd ones.
\end{theorem}

The values of \eqref{eq:thetaChar} are nothing but the branch points for the case of hyperelliptic curves, which directly deliver the image of the canonical map. In the case of non-hyperelliptic curves, \eqref{eq:thetaChar} gives the so-called \emph{multitangent hyperplanes} of the canonical model in $\mathbb{P}^{g-1}$. For instance, bitangent lines of smooth plane quartics in genus 3, tritangent planes of smooth space sextics in genus 4. It has been proven that the odd theta characteristics recover its algebraic curve~\cite{CapSer}. For explicit reconstructions the algebraic curve from their multitangent hyperplanes for small genera see~\cite{LehGenus4, LehGenus5, CelKulRenNamGenus4}.

Algebro-geometric solutions of integrable systems contribute solutions to the Torelli and the Schottky problems in any genus~\cite{Krichever1977, Dubrovin1981}, where fundamental objects are again the theta functions~\eqref{eq:thetaFunctionChar}. For our experiments, we use an implementation presented in~\cite{AgoCelEke} that follows these studies to recover curves from their Riemann matrices.

 Mathematical software packages are available to compute with theta functions, such as \cite{BruGan, AgoChu, DecEtAl, FraJabKle}, which enable us to carry out our experiments.

\subsection{Translation surfaces}\label{sec:Translation}

We will give two equivalent definitions of a translation surface, introduce the two families of examples considered in this paper, and conclude by discussing the connection to algebraic curves.

Fix $X$ a Riemann surface of genus $g$, and recall that $\omega$ is a \emph{holomorphic $1$-form} (also called an \emph{abelian differential}) if for every $x\in X$ there is a holomorphic function $f_x$ so that in local coordinates $\omega = f_x(z)\,dz$ with the condition that a transition map $T$ between charts with $f(z)$ and $g(z)$ as holomorphic function satisfies $f(T(z)) T'(z) = g(z)$. In other words $\omega$ is a global section of the cotangent bundle of $X$.
\begin{definition}\label{def:Translation}
    A \emph{translation surface} is:
    \begin{enumerate}
        \item a pairing $(X,\omega)$ where $X$ is a Riemann surface and $\omega$ is a holomorphic $1$-form. Two translation surfaces $(X,\omega)$ and $(Y,\eta)$ are equivalent if there exists a holomorphic diffeomorphism $\phi: X\to Y$  so that $\phi^*\eta= \omega$.
        
        \item collections of polygons up to an equivalence relation: $\mathcal{P}/\sim$. In particular $\mathcal{P}$ is a finite collection of polygons $\mathcal{P}$ in the plane so that all sides come in pairs of equal length with opposite orientations on the boundary of the polygons. Identifying these sides gives a compact finite genus Riemann surface. Given such collections of polygons $\mathcal{P}$ and $\mathcal{Q}$, we say $\mathcal{P}\sim \mathcal{Q}$ if elements of $\mathcal{P}$ can be cut into pieces along straight lines (where a cut produces two new boundary components that are paired) and these pieces can be translated and re-glued (where gluing only occurs along paired edges) to $\mathcal{Q}$.
    \end{enumerate}
\end{definition}
For more background on translation surfaces see \cite{Wright15, Massart22}.
\begin{remark}\label{rmk:zeroes}
    If $X$ has genus $g$ any holomophic $1$-form $\omega$ has $2g-2$ zeros with multiplicity. Away from the zeroes, $\omega = dz$, and if $z$ is a point where $\omega$ has a zero of order $k$, then we can write $\omega = z^k \,dz$. 
\end{remark}
\begin{example}
   The square torus $\mathcal{P} = [0,1]^2$ on the right and the torus $\mathcal{Q}$ on the left with vertices given by $(0,0), (1,0), (1,1), (2,1)$ are equivalent via the cut and paste operation shown below.
   \begin{center}
  \begin{tikzpicture}[line cap=round,line join=round,x=1.5 cm,y=1.5 cm]

\draw[-] (.5,0) -- (1.5,0)--(2.5,1)--(1.5,1)--(0.5,0);

\draw[-] (3,0) -- (4,0)--(5,1)--(4,1)--(3,0);
\draw[dashed] (4,1)-- (4,0);

\draw[-] (5.25,0) -- (6.25,0) -- (6.25,1)-- (5.25,0);
\draw[-] (7.5,1)--(6.5,1)-- (6.5,0)-- (7.5,1);

\draw[-] (8,0) --(9,0)--(9,1)--(8,1)--(8,0);
\draw[dashed] (8,0)--(9,1);
\begin{scriptsize}

\draw (1,0) node[anchor=north] {$1$};
\draw (2,1) node[anchor=south] {$1$};
\draw (1,.75) node[anchor=east] {$2$};
\draw (2, .25) node[anchor=west] {$2$};

\draw (3.5,0) node[anchor=north] {$1$};
\draw (4.5,1) node[anchor=south] {$1$};
\draw (3.5,.75) node[anchor=east] {$2$};
\draw (4.5, .25) node[anchor=west] {$2$};
\draw (4,.5) node[anchor = east]{$3$};

\draw (5.75,0) node[anchor=north] {$1$};
\draw (7,1) node[anchor=south] {$1$};
\draw (5.75,.75) node[anchor=east] {$2$};
\draw (7, .25) node[anchor=west] {$2$};
\draw (6.25,.5) node[anchor = east]{$3$};
\draw (6.5,.5) node[anchor = west]{$3$};

\draw (8.5,0) node[anchor=north] {$1$};
\draw (8.5,1) node[anchor=south] {$1$};
\draw (8.5,.55) node[anchor=east] {$2$};
\draw (8,.5) node[anchor = east]{$3$};
\draw (9,.5) node[anchor = west]{$3$};

\end{scriptsize}
\end{tikzpicture}
    \end{center}
\end{example}

The first definition of a translation surface is very concise, but the second definition is useful for constructing examples, and we will use these as our source of examples for this paper. As mentioned in Remark~\ref{rmk:zeroes}, there are $2g-2$ zeroes with multiplicity \cite[Theorem 1.2]{Wright15} If we label the zeroes by a multi-index $\alpha= (\alpha_1,\ldots, \alpha_k)$, then $k$ is the number of distinct zeroes, each with multiplicity $\alpha_k$, and $\sum_k \alpha_k= 2g-2$.
\begin{definition}
        The set of translation surfaces with orders of zeroes given by a multi-index $\alpha$ is the \emph{stratum $\mathcal{H}(\alpha)$}. When $k=2g-2$ and thus for each $1\leq i\leq 2g-2$, $\alpha_i = 1$, then the stratum is called the \emph{principal stratum.}
\end{definition}
\begin{example}
   Consider the $L$ shaped polygon of Figure~\ref{fig:LId}, which lives in the stratum $\mathcal{H}(2)$. 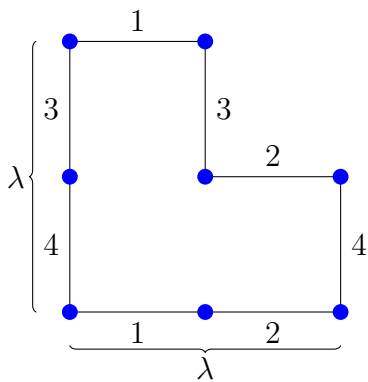
\begin{figure}[htbp]
    \centering
    \begin{tikzpicture}[line cap=round,line join=round,x=1.8cm,y=1.8cm]

\draw[-] (0,0) -- (2,0)--(2,1)--(1,1)--(1,2)--(0,2)--(0,0);
\draw (.5,0) node[anchor=north] {$1$};
\draw (.5,2) node[anchor=south] {$1$};
\draw (1.5,0) node[anchor=north] {$2$};
\draw (1.5,1) node[anchor=south] {$2$};
\draw (0,1.5) node[anchor=east] {$3$};
\draw (1,1.5) node[anchor=west] {$3$};
\draw (0,.5) node[anchor=east] {$4$};
\draw (2,.5) node[anchor=west] {$4$};

\fill  [color=blue] (0,0) circle (3pt);
\fill  [color=blue] (1,0) circle (3pt);
\fill  [color=blue] (1,1) circle (3pt);
\fill  [color=blue] (2,0) circle (3pt);
\fill  [color=blue] (1,2) circle (3pt);
\fill  [color=blue] (0,2) circle (3pt);
\fill  [color=blue] (0,1) circle (3pt);
\fill  [color=blue] (2,1) circle (3pt);

\draw [decorate,decoration = {brace}] (-.25,0) --  (-.25,2);
\draw (-.25,1) node[anchor=east] {$\lambda$};
\draw [decorate,decoration = {brace}] (2,-.25) --  (0,-.25);
\draw (1,-.25) node[anchor=north] {$\lambda$};
\end{tikzpicture}

    \caption{A symmetric $L$ translation surface with opposite sides identified by complex translation.  The vertices all map to the same point under identification. Using an Euler Characteristic argument, the underlying surface is a genus 2 Riemann surface. For a general symmetric $L$ shape fix $\lambda>1$ let $1$ be the length of sides $4$ and $1$, and let $\lambda - 1$ be the length of sides $2$ and $3$. The image above is shown for $\lambda = 2$}
    \label{fig:LId}
\end{figure}
   
   The sides given by the numbers $1,2,3,4$ are identified by complex translations $z\mapsto z \pm \mathbf{i}\lambda$, $z\mapsto z\pm \mathbf{i}$, $z\mapsto z \pm 1$, and $z\mapsto z\pm \lambda$, respectively. Under these edge identifications, the corners are all mapped to a single point, which is a zero with angle $6\pi = (1+2)(2\pi)$. An angle of $6\pi$ is a zero of order $2$ since it contains $2$ full circles of excess angle. This can also be seen in the local charts from $\mathbb{C}$ into the Riemann surface $z \mapsto z^{3}$, for which the differential is $3z^2 \,dz$, giving the zero of order $2$. Since there are $2g-2$ zeroes of $\omega$ with multiplicity, this implies that the genus should be $2$. We can also see this via the Euler characteristic $2-2g = V - E + F = 1 - 4 + 1$, where there is $1$ vertex under identification, 4 edges, and a single facet.
   
\end{example}

\begin{example}\label{ex:JSgenus2}
    The other primary examples we will work with are a family of curves in the principal stratum, so the $2g-2$ zeros are all distinct zeroes of order $1$, namely the \emph{Jenkins--Strebel differentials}, which are formed by a single horizontal rectangle of length $4g-4$ and height $1$ \cite{ZorichJS}. In Section~\ref{sec:JS}, we review these objects in more details. By considering horizontal lines connecting two zeros, there are at most $4g-4$ total parallel lines connecting all of the zeroes as explained in Remark~\ref{rem:4g-4zeros}. 
    %We will work with representatives called \emph{Jenkins--Strebel differentials} formed by a single horizontal rectangle of length $4g-4$ and height $1$ \cite{ZorichJS}. 
    In this example, the edge identifications are given through a permutation identifying the edges on the top of the rectangle to the bottom edges, as well as the horizontal translation identifying the vertical sides. For example in Figure~\ref{fig:g=2JS}, the translation surface $J_2$ is composed of a $1\times 4$ rectangle where the numbers indicate the sides identified by translation. There are $4$ horizontal sides on the top and bottom, and the sides are glued by taking a permutation of the top sides to glue to the bottom sides. By doing an Euler characteristic argument, we can verify that there are $2$ vertices, $5$ edges, and $1$ facet, resulting in a genus 2 surface. Thus $J_2 \in \mathcal{H}(1,1)$. Higher genus examples will be constructed in Section~\ref{sec:JS}
    
    By varying the edge lengths, we can move from $J_2$ through a family of curves which are all Jenkins--Strebel differentials of the same genus, similar to how $\lambda$ formed a one-parameter family of curves in $\mathcal{H}(2)$. In order to preserve some of the symmetries, we will only allow two parameters to change. Given $\lambda, \mu \in (0,\infty)$, let $0$ have side length $\lambda$, side $1$ has side length $\mu$, and sides $2,3$ and $4$ each have length $1$. This now gives a 2 dimensional family of curves $J_2(\lambda, \mu) \in \mathcal{H}(1,1).$ 
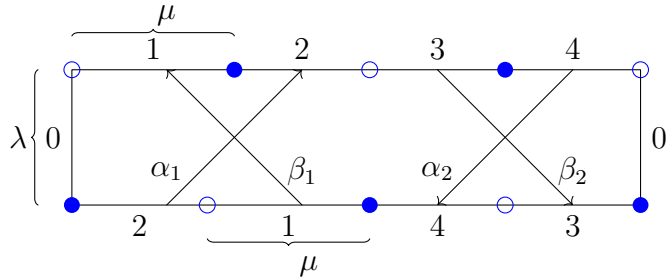
\begin{figure}[htbp]
    \centering
    \begin{tikzpicture}[line cap=round,line join=round,x=1.8cm,y=1.8cm]

\draw[-] (-.2,0) -- (2,0)--(4,0)--(4,1)--(2,1)--(-.2,1)--(-.2,0);
\draw (-.2,.5) node[anchor=east] {$0$};
\draw (4,.5) node[anchor=west] {$0$};

\draw (.4,1) node[anchor=south] {$1$};
\draw (1.5,1) node[anchor=south] {$2$};
\draw (2.5,1) node[anchor=south] {$3$};
\draw (3.5,1) node[anchor=south] {$4$};

\draw (.3,0) node[anchor=north] {$2$};
\draw (1.4,0) node[anchor=north] {$1$};
\draw (2.5,0) node[anchor=north] {$4$};
\draw (3.5,0) node[anchor=north] {$3$};
\fill  [color=blue] (-.2,0) circle (3pt);
\fill  [color=blue] (1,1) circle (3pt);
\fill  [color=blue] (2,0) circle (3pt);
\fill  [color=blue] (3,1) circle (3pt);
\fill  [color=blue] (4,0) circle (3pt);
\draw  [color=blue] (-.2,1) circle (3pt);
\draw  [color=blue] (.8,0) circle (3pt);
\draw  [color=blue] (2,1) circle (3pt);
\draw  [color=blue] (3,0) circle (3pt);
\draw  [color=blue] (4,1) circle (3pt);

\draw[->](.5,0)--(1.5,1);
\draw[color=black] (.5,.25) node {$\alpha_1$};
\draw[<-] (2.5,0)--(3.5,1);
\draw[color=black] (2.5,.25) node {$\alpha_2$};
\draw[->](1.5,0)--(.5,1);
\draw[color=black] (1.5,.25) node {$\beta_1$};
\draw[<-] (3.5,0)--(2.5,1);
\draw[color=black] (3.5,.25) node {$\beta_2$};

\draw [decorate,decoration = {brace}] (-.45,0) --  (-.45,1);
\draw (-.45,.5) node[anchor=east] {$\lambda$};
\draw [decorate,decoration = {brace}] (2,-.25) --  (.8,-.25);
\draw (1.55,-.3) node[anchor=north] {$\mu$};
\draw [decorate,decoration = {brace}] (-.2,1.25) --  (1,1.25);
\draw (.5, 1.25) node[anchor=south] {$ \mu$};
%\draw[color=black] (-0.39,-0.22) node {$n$};

\end{tikzpicture}
    \caption{Above is $J_2(\lambda,\mu)$ given by the $4$ unit squares glued with the following identifications. The associated permutation is $\pi_{J_2} = \begin{pmatrix}1&2&3&4\\ 2 & 1 & 4&3\end{pmatrix}$. The paths $\alpha_j, \beta_j$ for $j=1,2$ form a symplectic basis of homology.}
    \label{fig:g=2JS}
\end{figure}
\end{example}

%The benefit of translation surfaces from the perspective of this project is they provide a way of systematically constructing examples of Riemann surfaces. In order to understand a Riemann surface as an algebraic curve, we must have some basis of one-forms $\omega_1,\ldots, \omega_g$. The difficulty is that these polygonal constructions automatically come with a holomorphic $1$-form inherited from $dz$, making it difficult to find the other $g-1$ basis elements for the underlying Riemann surface.

For the Jenkins--Strebel differentials, we do not know the underlying Riemann surface. However in the case of the $L$, we can use symmetries of the polygon to determine a second linearly independent $1$-form. This was done by \cite{Sil, Rod} in the case of the $L$ using the order $4$ symmetry to show that the underlying Riemann matrix is given by
\begin{equation}\label{SilholCurve}
    \tau_{\lambda}= 
\frac{\mathbf{i}}{2\lambda - 1}\begin{pmatrix}2\lambda^2 - 2\lambda + 1 & -2\lambda(\lambda-1) \\ -2\lambda(\lambda-1) & 2\lambda^2 -2\lambda + 1\end{pmatrix}.
\end{equation}
Moreover, since we are in genus 2 and thus hyperelliptic, the equation of the underlying curve is given by
$$y^2= x(x^2- 1)(x-a)(x-1/a) \text{ for } a\neq -1,0,1.$$
The values of $a$ are computed in \cite{Sil} for certain values of $\lambda$, in particular for the examples computed in this paper that when $\lambda = 2$, then $a= 7 + 4\sqrt{3}$. Other families of surfaces where the symmetries are used to compute the underlying curve can be found in \cite{DM21}.

 We close our section with an example that illustrates reconstructing an algebraic curve from its translation surface via numerical computations relying on methods in Section~\ref{sec:RS}.

 \begin{example}\label{ex:SilholLambda2} We take the Riemann matrix $\tau$ as $\mathbf{i} \begin{pmatrix} \frac{5}{3} & -\frac{4}{3} \\ -\frac{4}{3} & \frac{5}{3} \end{pmatrix}$ for certain $\lambda$ in Equation~\eqref{SilholCurve}. Section~\ref{sec:Translation} contains more details about this example and its underlying translation surface. We use SageMath~\cite{BruGan}, computing with 100 bits of precision, we approximate the six branch points as follows via the six theta constants with odd characteristics: 
 
 \begin{center}
      \begin{small}
 $\begin{matrix}
 -2.0000000000000000000000000000 + \mathbf{i}6.4112869792140406597766726185\cdot 10^{-62},\\
 -1.0242537764555949265655782388 + \mathbf{i}4.6398778498086499909043081781 \cdot 10^{-64},\\
 -0.50000000000000000000000000001 - \mathbf{i}3.0730068477907671248074738783 \cdot 10^{-62},\\
 -0.97632053987682181152178995357 - \mathbf{i}3.8173070449992676299965669304 \cdot 10^{-64},\\
 1.2417360351295279957623671524 + \mathbf{i}1.2694991939112402714970349296 \cdot 10^{-61},\\
 0.80532413629736369749980162893 + \mathbf{i}1.1577740944972796938745848090\cdot 10^{-62}.
 \end{matrix}$
 \end{small}
  \end{center}
In fact, the branch points are nothing but the quotients of the six theta constants $-\theta^{\varepsilon,\delta}_{1}/\theta^{\varepsilon,\delta}_{2}$ among all the even characteristics with the notation~\eqref{thetaConstant}. To verify that this curve is isomorphic to the curve exhibited in~\eqref{SilholCurve} one can compute their absolute Igusa invariants in a SageMath class~\cite{SageIgusa} and note that the absolute Igusa invariants of both curves coincide up to a numerical round-off.
\end{example}

\subsection{Discrete Riemann surfaces and discrete period matrices}\label{sec:DiscreteApproximations}

The section aims to assist the reader in the literature on discrete Riemann surfaces. Our references are \cite{BobGun, Mercat2001}. We here skip recalling the vast amount of technical background on the topic and pinpoint the related results in the references instead. 

Given a surface $X$, namely a two real-dimensional manifold, one discretizes the surface by considering it via one of its cellular decompositions, say $\Lambda$ together with a discrete complex structure. Here, the discrete complex structure is introduced with the consideration of the dual cell decomposition of $\Lambda$, denoted by $\Lambda^*$. Note that \cite{BobGun} encodes these objects, namely $\Lambda$ and $\Lambda^*$, by the colors black and white, respectively. Set $\Gamma:=\Lambda \cup \Lambda^*$, which is called the double of the cell decomposition $\Lambda$. This is to set a theoretical framework for \emph{discrete complex analysis}, that is to say the discrete theory of complex analytic functions. 

First and foremost, the discretization of the Cauchy-Riemann equation is formulated in terms of combinatorial elements of $\Lambda$ and its dual $\Lambda^*$, namely the sub-cells. In particular, for a complex valued function $f$ defined on the $0$-cells of $\Gamma$ to satisfy the Cauchy-Riemann equation means certain compatibility between the proportions of the 0-cells of $\Lambda$ and $\Lambda^*$ and their values under $f$, see~\cite[Section 3]{BobGun}.

In the development of this theory, one fundamental concept is the discrete theory of Riemann surfaces. In fact, many results in the classical theory have discrete counterparts, which includes period matrices, Abelian integrals, and so forth. Bobenko and Günther use a medial graph approach to the discrete theory of Riemann surfaces on quad-decompositions, which was introduced as a perspective to discrete complex analysis in \cite{Bobenko2016}. Mercat makes use of the tool of deRham cohomology to introduce standard notions in discrete exterior calculus~\cite{Mercat2001}. A \emph{discrete one-form} $\omega$ is defined as a complex function on the one-cells of $\Gamma$. The evaluation of $\omega$ at an oriented edge is nothing but the discrete integral $\int_e\omega$, via which one defines the \emph{discrete integral} over a directed path by means of the oriented edges forming the path. 

Suppose that $X$ is a compact Riemann surface of genus $g$. Its discretization is given a cell decomposition of $X$. We follow~\cite{BobGun} for the notion of discrete Riemann matrix, where the authors note that their object coincides with~\cite{Mercat2001}. Fix a symplectic basis $\alpha_1,\dots,\alpha_g,\beta_1,\dots, \beta_g$ of ${\rm H}_1(X,\mathbb{Z})$. One defines the $A$ and $B$ periods of a given discrete differential by taking the integrals over the $\alpha_j$ and $\beta_j$, respectively. These cycles induce closed paths on $\Lambda$ and $\Lambda^*$, which are distinguished with the colors, white and black, labeling the cell decompositions. This yields the notions of black or white $A$ or $B$ periods. 

For a discrete differential, the $4g$ discrete black and white \emph{periods} are defined as the integrals that are over the induced black and white closed paths. For technical details of the cycles and the periods, see~\cite[Section 5.1]{BobGun}. It turns out there is a unique holomorphic differential such that the black and white $A$ and $B$ periods match a given set of $4g$ complex values \cite[Theorem 6.3]{BobGun}. It follows that the \emph{canonical basis} of holomorphic one-forms $\omega_1,\ldots, \omega_g$ \cite[\S 6.3 Definition]{BobGun} is well defined where the black and white $A$ periods are chosen to be equal with integration against the curves $\alpha_1,\ldots, \alpha_g$ is the identity matrix. The $g\times g$ \emph{discrete period matrix} entries are the $B$ periods with respect to the canonical basis. This is somewhat mimicking the normalization of the $g\times 2g$ period matrix in the classical setting. Abusing the notation, we will call the discrete period matrix as \emph{discrete Riemann matrix} by referring its second $g\times g$ part.  

There is another notion of a period matrix called the \emph{complete discrete period matrix}, which is a $2g\times 2g$ block matrix made of four $g\times g$ matrices. These $g\times g$ matrices are formed by not imposing the white and black $A$ periods are equal, and instead considering $g\times g$ matrices formed according to the relationship of black and white periods. Note that the discrete period matrix can be computed from the complete one~\cite[Remark at Page 917]{BobGun}. 

%Given a discrete Riemann surface, that is a discretization of the underlying Riemann surface with a cell decomposition, 
Computing the (complete) discrete period matrix amounts to computing the discrete periods. The periods are the values of the discrete differentials at the edges of the closed paths arising from the fixed symplectic basis. In order to compute these values, one may use the condition of being holomorphic for the discrete differentials by the discrete Cauchy-Riemann relations. This gives a linear system of equations, which we call \emph{holomorphicity equations}. 

In our algorithm we also construct the so-called \emph{periodicity equations}, which are given by the presentation of our underlying surface as a translation surface. Translation surfaces are an example of a \emph{polyhedral surface}, which consists of planar polygons that are glued together along edges. It turns out that this is yet another characterization of a compact Riemann surface~\cite{Bost}. This perspective might be more convenient for explicit computations involving discrete surfaces, in particular while considering the discrete complex structure on the decomposition of the surface. For the case of computing the discrete period matrices, the edges being glued adds linear equations to the holomoprhicity equations, namely the periodicity equations.  

As the decomposition into cells gets finer, one expects that the discrete period matrix converges to a Riemann matrix of the underlying Riemann surface. At this point the convergence of discrete Riemann surfaces to their continuous counterparts in full generality remains open. When one decomposes the surface into a Delaunay triangulation, convergence of the period matrices is proved.
\begin{theorem}[Theorem 2.5, \cite{BobSko}] For a sequence of triangulations of a compact Riemann surface $X$ with the maximal edge length tending to zero and with face angles bounded from zero, the discrete period matrices converge to a Riemann matrix of $X$. 
\end{theorem}
The techniques of \cite{BobSko} led to the suggestion that one builds the discrete complex structure via quad-decompositions with orthogonal diagonals. \cite{BobGun} put the suggestion in place and examined related fundamental notions of such discrete Riemann surfaces. Another convergence result of~\cite[Theorem 3]{BobBuc} considers compact Riemann surfaces in terms of their branched covers of the Riemann spheres. Our experiments, which employ compact Riemann surfaces given by quadrangulations, give evidence for a potential result on convergence of the period matrices of \cite{BobGun}.

The convergence rate is given by looking at particular triangulations as the longest edge length $l$ goes to zero. For translation surfaces, the convergence rate \cite[Theorem 2.5]{BobSko}, is proportional to $l$ when the genus $g=1$. When $g\geq 2$, the rate is $l|\log(l)|$ when the $2g-2$ zeroes are all distinct, and at worst $l^\frac{2}{2g-1}$ which corresponds to the case when there is a single zero of order $2g-2$. In \cite[Theorem 2.4]{BobBuc}, when working with ramified coverings of the sphere and carefully choosing an optimal triangulation, the convergence rate is always proportional to $l$. We consider square quadrangulations instead of triangulations for two primary reasons. First, not every surface can be split into quadrilaterals, but when this is possible, the computational ease of parameterizing the quadrangulation, subdivision of the quadrilaterals, and elementary linear algebra needed to compute the matrices makes it more effective. Moreover, after discussions with Felix G{\"u}nther, we expect similar convergence rates to \cite{BobSko}. The key idea is that by dividing the squares along the diagonal, the resulting triangulation is sufficiently regular to fit into the framework of \cite{BobSko}. 

\section{Algorithms} \label{sec:algorithms}

In the following section, we will present algorithms for constructing discrete Riemann matrices associated to two families of translation surfaces: The $L$ shape in Section~\ref{sec:Lsection} and the Jenkins--Strebel representatives in Section~\ref{sec:JS}. In the first case, we aim to construct the associated discrete Riemann matrices when approximating the shape by squares and observe the convergence to the known underlying algebraic curve. In the second case we observe convergence and show experiments indicating what we expect from the underlying curve.

\begin{definition}\label{def:levels}
        Given a translation surface $P$, the \emph{$0$th level discrete approximation}  $P_0$ is the discrete Riemann surface defined by the smallest bipartite square quadrangulation of the surface which respects the identifications of the vertices in $P$. The $n$th level discrete approximation $P_n$ for $n\in \mathbb{N}_{\geq 1}$ is defined to be the discrete Riemann surface defined by the bipartite square quadrangulation which subdivides each square of the $n=0$ level into $3^{2n}$ squares.
    \end{definition}
    
    We give an algorithm used to obtain the discrete Riemann matrix associated to $P_n$.   
\subsection{General Polygon Algorithm} \label{alg:general}
\noindent \textbf{Input:} A polygon $P$ in $\mathbb{R}^2$ with fixed side lengths, and side identifications which gives a translation surface, and a level $n\in \mathbb{N}_{\geq 0}$.

\noindent \textbf{Output:} The discrete Riemann matrix associated to $P_n$.

\begin{enumerate}
    \item \textbf{Constructing an initial bipartite quadrangulation.}  We first divide the given translation surface into a bipartite quadrangulation, note that the bipartite quadrangulation must also respect the identifications of the vertices under the identifications in the given translation surface.
    \item \textbf{Quadrangulations for further levels of approximation.} In order to preserve the bi-coloring, we must divide each square into an odd number of squares, for which we will choose $3$. So each square of side length $s$ will be divided into $3^{2n}$ squares, and so the $n$th level approximation will consist of squares of size $ {s}/{3^n}.$
    
    \item \textbf{Labelling vertices.} Placing the bottom left corner of the given polygon at the origin, we label the the vertices by their location in the plane, where we will move freely between complex notation and vector notation:
    $$x_{i,j}= \left(i\frac{s}{3^n }, j\frac{s}{3^n } \right) = i\frac{s}{3^n} +  j \frac{s}{3^n} \mathbf{i} .$$ The vertex $x_{i, j}$ corresponds to the bottom left corner of the square which is $(i+1)$ from the left and $(j+1)$ from the bottom.
    
     \item \textbf{Holomorphicity equations.} For each square with coordinate $x_{i,j}$ in the bottom left, we have the holomorphicity relation
    \begin{equation}\label{eq:holo}
        \mathbf{i}(x_{i+1,j+1} - x_{i,j}) = x_{i,j+1} - x_{i+1,j}.
    \end{equation}
     See Figure \ref{fig:holo} for a visual description of these relations.
    
    \begin{figure}[htbp]
    \centering
    \begin{tikzpicture}[line cap=round,line join=round,x=1.5cm,y=1.5cm]

\draw[-] (0,0) -- (2,0)--(2,2)--(0,2)--(0,0);
\draw (0,0) node[anchor=north east] {$x_{i,j}$};
\draw (2,0) node[anchor=north west] {$x_{i+1,j}$};
\draw (0,2) node[anchor=south east] {$x_{i,j+1}$};
\draw (2,2) node[anchor=south west] {$x_{i+1,j+1}$};
\fill  [color=blue] (2,2) circle (3pt);
\fill  [color=blue] (0,0) circle (3pt);
\draw  [color=blue] (0,2) circle (3pt);
\draw  [color=blue] (2,0) circle (3pt);

\draw[->](.25,1.75)--(1.75,.25);
\draw[->] (1.75,1.75)--(.25,.25);

%\draw[color=black] (-0.39,-0.22) node {$n$};

\end{tikzpicture}
    \caption{Representation of the holomorphicity equations at the $ij$ square.}
    \label{fig:holo}
\end{figure}
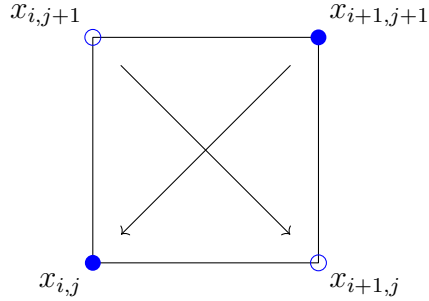
By \cite[Proposition~2.1]{BobGun} and the discrete Cauchy Riemann equations, a discrete Riemann surface tiled by squares with given orientation in Figure~\ref{fig:holo} should indeed satisfy \eqref{eq:holo}. When $x_{i,j}$ is instead a white vertex \cite{BobGun} gives
$$\mathbf{i}  = \frac{x_{i+1,j+1}- x_{i,j}}{x_{i+1,j} - x_{i,j+1}} =\frac{x_{i+1,j+1}- x_{i,j}}{-(x_{i,j+1} - x_{i+1,j})} $$
which is equivalent to~\eqref{eq:thetaFunctionChar}.

  \item \textbf{Periodicity.} We choose a symplectic basis of homology $\alpha_1,\ldots, \alpha_g,\beta_1,\ldots, \beta_g$, and the associated $4g$ discrete periods are given by $A_k^w, A_k^b, B_k^w, B_k^b$ for $k=1,\ldots, g$, and the superscripts $w,b$ representing the white and black periods, respectively. The periodicity relationships are constructed in order to make edge identifications for the given polygon in the plane. So for each edge identification, the identified vertices have difference given by the correct associated vector, and the coloring is found through checking the parity of $i$ and $j$.

    \item\textbf{Final Normalizations.} To make a well-determined system, we make the following normalizations:
    \begin{itemize}
        \item Fix the first values of the holomorphic function, one one black and one one white vertex: $x_{0,0} = x_{1,0} = 0$. This comes from the fact that the associated holomorphic function is only defined up to a constant, so we normalize the constant to be zero at the origin.
        \item In order to construct the canonical basis of discrete holomorphic differentials, we set the black and white values to be the same, and for each $k=1,\ldots, g$, the differential $\omega_k$ is determined by the following equations for $j=1,\ldots, g$: $$A_j^w = A_j^b \quad \text{and}\quad A_{j}^w = \begin{cases} 1 & j = k\\
        0 & \text{else.}\end{cases}$$
    \end{itemize}
    
    \item \textbf{Solving a system of equations for the discrete approximation.} Now for each $k=1,\ldots, g$, the $k$th row of the discrete period matrix is given by
    $$\frac{1}{2}(B_1^w + B_1^b,\ldots, B_g^w + B_g^b).$$
\end{enumerate}

\subsection{Riemann Matrix of the \texorpdfstring{$L$}{L}} \label{sec:Lsection}

In this section we find numerical approximations to a family of translation surfaces for which we already know the Riemann matrix. Namely, consider a symmetric $L$ shape with side length $\lambda \in (1,\infty)$ with opposite sides identified as in Figure~\ref{fig:LId}. We aim to construct the associated discrete Riemann matrices when approximating the shape by squares, and observe the convergence towards $\tau_\lambda$ in Equation~\eqref{SilholCurve}, which is guaranteed by~\cite{BobGun}.

%\subsubsection{An algorithm for rational $\lambda$}

As mentioned in Section~\ref{sec:DiscreteApproximations}, to construct the discrete Riemann matrix, we must solve a system of linear equations given by holomorphicity relations between vertices of squares in the polygonal subdivision, and by periodicity relations obtained from identifications of points on the boundary of the shape. We now describe an algorithm which constructs for us this system of equations. See~\cite{BobGun} for more details on this construction. We used MATLAB for the implementation of the algorithm.

\subsubsection{Algorithm for symmetric \texorpdfstring{$L$}{L}} \label{alg:L}
\noindent\textbf{Input:} Let $\lambda = {p}/{q}$ be rational and reduced so that $\gcd(p,q) = 1$. Let $n \in \mathbb{N}\cup{\{0\}}$ be the level of approximation. 

\noindent\textbf{Output:} Discrete Riemann matrix of the $n$th level approximation for symmetric $L$ with side length $\lambda$.

\begin{enumerate}
    \item \textbf{Constructing an initial bipartite quadrangulation.}  Refer to Figure \ref{fig:labels} for an example of the level $0$ approximation. To divide the entire shape into squares, the sizes of the squares must divide ${1}/{q}$. In order to bi-color the square tiling and maintain the vertex identifications, there must be an even number of squares on each side length. Hence we define the step size to be $ s_\lambda = \mathrm{lcm}(q,2).$ The shape $L$ can be bi-colored by being divided into squares of size $1/s_\lambda$.
    \item \textbf{Quadrangulations for further levels of approximation.} Each square of side length ${1}/{s_\lambda}$ will be divided into $3^{2n}$ squares, and so the $n$th level approximation will consist of squares of size $ {1}/({3^ns_{\lambda}}).$ 
    \item \textbf{Labelling vertices.} We label the vertices with the following bounds to match the $L$ shape:
    $$x_{i,j}= \left(\frac{i}{3^n s_\lambda}, \frac{j}{3^n s_\lambda} \right) \quad \text{ for \quad } \begin{cases} 0\leq i \leq  3^n  s_\lambda & 0\leq j \leq \lambda 3^n s_\lambda, \\
    3^n s_\lambda +1 \leq i \leq \lambda 3^n s_\lambda & 0 \leq j \leq 3^n s_\lambda.
    \end{cases}$$
    
    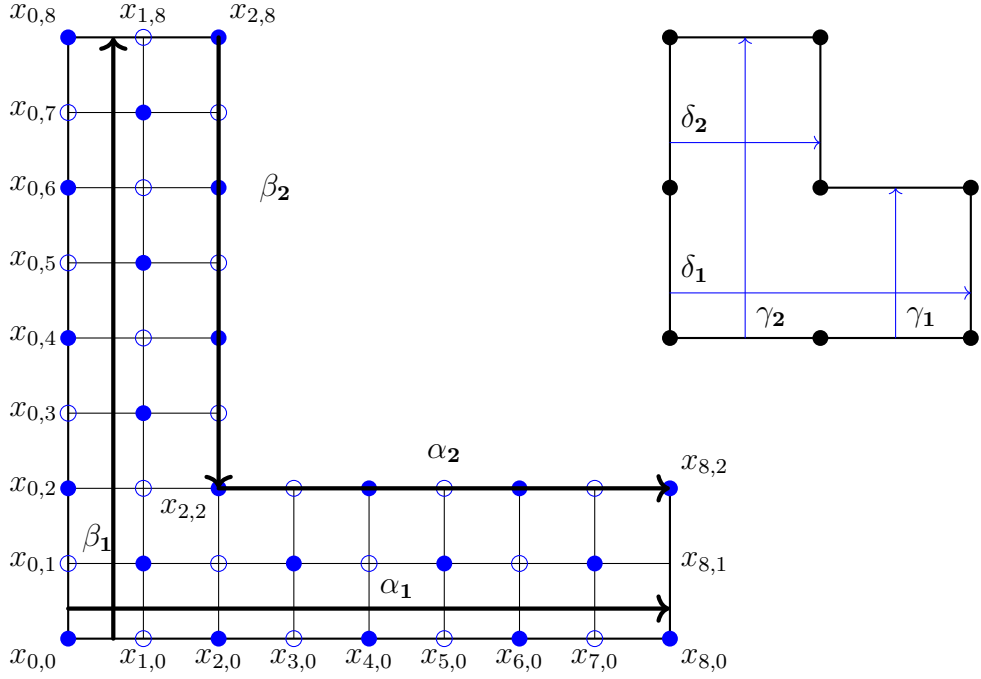
\begin{figure}[htb]
    \centering
    \begin{tikzpicture}[line cap=round,line join=round,x=2 cm,y=2 cm]

\draw[-,thick] (0,0) --(4,0)--(4,1)--(1,1)--(1,4)--(0,4)--(0,0);
\draw[-] (1,0) -- (1,1)-- (0,1);
\draw[-] (.5,0) -- (.5,4);
\draw[-] (0,.5) -- (4,.5);
\draw[-] (1.5,0) -- (1.5,1);
\draw[-] (2,0) -- (2,1);
\draw[-] (2.5,0) -- (2.5,1);
\draw[-] (3.5,0) -- (3.5,1);
\draw[-] (3,0) -- (3,1);
\draw[-] (0,3.5) -- (1,3.5);
\draw[-] (0,3) -- (1,3);
\draw[-] (0,2.5) -- (1,2.5);
\draw[-] (0,2) -- (1,2);
\draw[-] (0,1.5) -- (1,1.5);

\draw (0,0) node[anchor=north east] {$x_{0,0}$};
\draw (.5,0) node[anchor=north ] {$x_{1,0}$};
\draw (1,0) node[anchor=north ] {$x_{2,0}$};
\draw (1.5,0) node[anchor=north ] {$x_{3,0}$};
\draw (2,0) node[anchor=north ] {$x_{4,0}$};
\draw (2.5,0) node[anchor=north ] {$x_{5,0}$};
\draw (3,0) node[anchor=north ] {$x_{6,0}$};
\draw (3.5,0) node[anchor=north ] {$x_{7,0}$};
\draw (4,0) node[anchor=north west] {$x_{8,0}$};
\draw (0,.5) node[anchor=east ] {$x_{0,1}$};
\draw (0,1) node[anchor=east ] {$x_{0,2}$};
\draw (0,1.5) node[anchor=east ] {$x_{0,3}$};
\draw (0,2) node[anchor=east ] {$x_{0,4}$};
\draw (0,2.5) node[anchor=east ] {$x_{0,5}$};
\draw (0,3) node[anchor=east ] {$x_{0,6}$};
\draw (0,3.5) node[anchor=east ] {$x_{0,7}$};
\draw (0,4) node[anchor=south east ] {$x_{0,8}$};

\draw (.5,4) node[anchor=south ] {$x_{1,8}$};
\draw (1,4) node[anchor=south west] {$x_{2,8}$};
\draw (4,.5) node[anchor=west ] {$x_{8,1}$};
\draw (4,1) node[anchor=south west] {$x_{8,2}$};
\draw (1,1) node[anchor=north east] {$x_{2,2}$};

\fill  [color=blue] (0,0) circle (3pt);
\fill  [color=blue] (1,0) circle (3pt);
\fill  [color=blue] (2,0) circle (3pt);
\fill  [color=blue] (3,0) circle (3pt);
\fill  [color=blue] (4,0) circle (3pt);
\fill  [color=blue] (0,1) circle (3pt);
\fill  [color=blue] (1,1) circle (3pt);
\fill  [color=blue] (2,1) circle (3pt);
\fill  [color=blue] (3,1) circle (3pt);
\fill  [color=blue] (4,1) circle (3pt);
\fill  [color=blue] (0,2) circle (3pt);
\fill  [color=blue] (1,2) circle (3pt);
\fill  [color=blue] (0,3) circle (3pt);
\fill  [color=blue] (1,3) circle (3pt);
\fill  [color=blue] (0,4) circle (3pt);
\fill  [color=blue] (1,4) circle (3pt);
\fill  [color=blue] (.5,.5) circle (3pt);
\fill  [color=blue] (1.5,.5) circle (3pt);
\fill  [color=blue] (2.5,.5) circle (3pt);
\fill  [color=blue] (3.5,.5) circle (3pt);
\fill  [color=blue] (.5,1.5) circle (3pt);
\fill  [color=blue] (.5,2.5) circle (3pt);
\fill  [color=blue] (.5,3.5) circle (3pt);

\draw  [color=blue] (0,.5) circle (3pt);
\draw  [color=blue] (0,1.5) circle (3pt);
\draw  [color=blue] (0,2.5) circle (3pt);
\draw  [color=blue] (0,3.5) circle (3pt);
\draw  [color=blue] (.5,0) circle (3pt);
\draw  [color=blue] (.5,1) circle (3pt);
\draw  [color=blue] (.5,2) circle (3pt);
\draw  [color=blue] (.5,3) circle (3pt);
\draw  [color=blue] (.5,4) circle (3pt);
\draw  [color=blue] (1,.5) circle (3pt);
\draw  [color=blue] (1,1.5) circle (3pt);
\draw  [color=blue] (1,2.5) circle (3pt);
\draw  [color=blue] (1,3.5) circle (3pt);
\draw  [color=blue] (1.5,0) circle (3pt);
\draw  [color=blue] (1.5,1) circle (3pt);
\draw  [color=blue] (2,.5) circle (3pt);
\draw  [color=blue] (2.5,0) circle (3pt);
\draw  [color=blue] (2.5,1) circle (3pt);
\draw  [color=blue] (3,.5) circle (3pt);
\draw  [color=blue] (3.5,0) circle (3pt);
\draw  [color=blue] (3.5,1) circle (3pt);

%\draw[->, ultra thick] (.2,0) --(.2, 4) ;
%\draw (.2,2) node[anchor= south west] {$\mathbf{\beta_1}$};
\draw[->, ultra thick] (0.3,0) --(0.3, 4) ;
\draw (0.02,.5) node[anchor= south west] {$\mathbf{\beta_1}$};
%\draw[->, ultra thick] (1,4) --(1, 1) ;
%\draw (1.25,3) node[anchor= west] {$\mathbf{\beta_2}$};
\draw[->, ultra thick] (1,4) --(1, 1) ;
\draw (1.2,3) node[anchor= west] {$\mathbf{\beta_2}$};
%\draw[->, ultra thick] (1,1) --(4, 1) ;
%\draw (2.5,1) node[anchor= south] {$\mathbf{\alpha_2}$};
\draw[->, ultra thick] (1,1) --(4, 1) ;
\draw (2.5,1.1) node[anchor= south] {$\mathbf{\alpha_2}$};

\draw[->, ultra thick] (0,.2) --(4,.2) ;
\draw (2,.2) node[anchor= south west] {$\mathbf{\alpha_1}$};

\draw[-,thick] (4,2) --(6,2)--(6,3)--(5,3)--(5,4)--(4,4)--(4,2);
\fill  [color=black] (4,2) circle (3pt);
\fill   [color=black] (6,2) circle (3pt);
\fill   [color=black] (6,3) circle (3pt);
\fill   [color=black] (5,3) circle (3pt);
\fill   [color=black] (5,4) circle (3pt);
\fill   [color=black] (4,4) circle (3pt);
\fill   [color=black] (5,2) circle (3pt);
\fill   [color=black] (4,3) circle (3pt);

\draw[->, color= blue] (4,2.3) --(6,2.3) ;
\draw (4,2.3) node[anchor= south west] {$\mathbf{\delta_1}$};
\draw[->, color= blue] (4,3.3) --(5,3.3) ;
\draw (4,3.3) node[anchor= south west] {$\mathbf{\delta_2}$};
\draw[->, color= blue] (4.5,2) --(4.5,4) ;
\draw (4.5,2) node[anchor= south west] {$\mathbf{\gamma_2}$};
\draw[->, color= blue] (5.5,2) --(5.5,3) ;
\draw (5.5,2) node[anchor= south west] {$\mathbf{\gamma_1}$};
\end{tikzpicture}
    \caption{This shows some of the labels in the level 0 approximation of a symmetric $L$ translation surface with opposite sides identified by complex translation.  In this case, we have $\lambda=4$. On the inset to the right, we show a basis of homology $(\gamma_1, \gamma_2, \delta_1, \delta_2)$, which gives a symplectic basis of homology $( \gamma_1, \gamma_2 - \gamma_1, \delta_1, \delta_2)$. In the left picture is the basis of homology that we used in our algorithm (see Remark \ref{rem:basis}).}
    \label{fig:labels}
\end{figure}

    \item \textbf{Holomorphicity equations.} For each bottom left of a square, we have a new holomorphicity equation. So in this case, we have a total of $3^{2n} (s_\lambda)^2 (2\lambda   - 1)$ equations with indices given by
     $$ \begin{cases} 0\leq i \leq  3^n s_\lambda - 1 & 0\leq j \leq \lambda 3^n s_\lambda -1, \\
    3^n s_\lambda  \leq i \leq \lambda 3^n s_\lambda -1 & 0 \leq j \leq 3^n s_\lambda - 1.
    \end{cases}$$
  
    \item \textbf{Periodicity equations.} We first choose a symplectic basis of the underlying Riemann surface, as shown in Figure~\ref{fig:labels}. To justify our choice, under the edge identifications, we select the closed loops $\alpha_1$ and $\beta_1$ as the first two, and normalize the symplectic basis so that we travel from $\alpha_1$ counterclockwise to $\beta_1$ for a positive intersection number. The choice of $\alpha_2$ comes naturally by trying to find another curve parallel to $\alpha_1$ which does not intersect $\beta_1$. The final step involves finding $\beta_2$. To do this, recall $\alpha_2$ is also identified by travelling from $x_{2,0}$ to $x_{8,0}$. Now going counterclockwise from $\alpha_2$, we travel around the vertex $x_{2,0}$, which is identified with $x_{2,8}$, so $\beta_2$ must travel from $x_{2,8}$to $x_{2,2}$ in order to have the correct intersection number with $\alpha_2$ and avoid intersecting $\beta_1$ or $\alpha_1$.
    
    We now need to construct the period equations, so for example to travel from $x_{0,0}$ to $x_{8,0}$, we travel along the $\alpha_1$ curve on black periods
    $$A_1^b = \int_{\alpha_1}\omega_k = x_{8,0} - x_{0,0}.$$
    
    For all the $A$ periods we have the following equations where the parity $p$ is determined by $p = b$ if $i+j \equiv 0 \pmod 2$ and $p=w$ otherwise:
    $$\begin{cases} 0\leq j \leq 3^n s_\lambda  & x_{\lambda 3^ns_{\lambda},j} - x_{0,j} = A_1^{p}\\
    3^{n} s_{\lambda} \leq j\leq \lambda 3^n s_{\lambda} & x_{3^n s_{\lambda},j} - x_{0,j} =  A_1^p - A_2^p .\end{cases}$$
    
   We compute similar equations for the $B$ periods, 
   $$\begin{cases} 0\leq i \leq 3^n s_\lambda  & x_{i,\lambda 3^ns_{\lambda}} - x_{i,0} = B_1^{p} \\
    3^{n} s_{\lambda} \leq i\leq \lambda 3^n s_{\lambda} & x_{i, 3^n s_{\lambda}} - x_{i,0} = B_1^p + B_2^p .\end{cases}$$
    In total we have $ 2(\lambda 3^n s_\lambda + 1)$ equations.

    \item\textbf{Final Normalizations.} In the final two normalizations we have the following number of equations:
    \begin{itemize}
        \item 2 equations for normalization of holomorphic function.
        \item For each $k=1,2$, there are $4$ equations normalizing to the canonical basis. 
    \end{itemize}
    
    \item \textbf{Solving a system of equations for the discrete approximation.} 
    For the $k$th row of the period matrix with $k=1,2$, we obtain the equations by solving the system with:
    \begin{itemize}
        \item Total of $9+ 3^{2n}(s_\lambda)^2 (2\lambda-1) + 2\lambda 3^n s_\lambda$ variables. With $(3^ns_\lambda +1)(\lambda 3^n s_\lambda +1) + (\lambda 3^n s_\lambda - 3^ns_\lambda) (3^n s_\lambda + 1)$ variables $x_{i,j}$, and $8$ variables coming from $B_j^p$ with $j=1,2$ and parity given by $b$ and $w$. 
        \item Total of $10 + 3^{2n}(s_\lambda)^2(2\lambda-1) + 2\lambda 3^n s_\lambda$ equations. With $3^{2n}(s_\lambda)^2(2\lambda-1)$ holomorphicity equations, $2(\lambda 3^n s_\lambda + 2)$ periodicity equations, and $6$ normalizing equations.
    \end{itemize}
    Thus we have a system of equations overdetermined by 1 equation, and these are not conflicting with a unique solution, as guaranteed by \cite[Theorem~6.8]{ BobGun} since they are simply relations describing the unique holomorphic differential with the given initial conditions.
\end{enumerate}

\begin{remark}\label{rem:basis}
    In our implementation, we choose the homology basis that is exhibited in Figure~\ref{fig:labels}, namely $\alpha_1,\alpha_2,\beta_1,\beta_2$. This basis appears to be in relative homology, but this just reflects the representation of the $L$ as a translation surface which includes a holomorphic one-form. Alternatively, in the same figure, we point out a different symplectic basis for the object given by $\delta_1,\delta_2,\gamma_1, \gamma_2 - \gamma_1$. The algorithm can be designed using the new basis by replacing the cases of the periodicity conditions as follows: 
    \begin{align*}
        x_{3^n s_{\lambda},j} - x_{0,j} &=   A_2^p, \\
        x_{i,\lambda 3^ns_{\lambda}} - x_{i,0} &= B_2^p + B_1^p,\\
        x_{i, 3^n s_{\lambda}} - x_{i,0} &= B_1^p.
    \end{align*}
    The matrix which results from $\delta_1,\delta_2,\gamma_1, \gamma_2 - \gamma_1$ can be obtained from the matrix resulting from $\alpha_1,\alpha_2,\beta_1,\beta_2$ by the transformation $$\begin{bmatrix}
        x & y\\ y & z
    \end{bmatrix} \mapsto \begin{bmatrix}
        x+2y+z & -y-z\\ -y-z & z
    \end{bmatrix}.$$
    For each of the two the symplectic bases, we approximated the corresponding Riemann matrices and their curves for several examples. We compared the results via their invariants, which coincide up to 15 digits. Though the same curve is represented under this change of basis, we note that the experiments arising from $\alpha_1,\alpha_2,\beta_1,\beta_2$ coincide with the results following the algebraic computations of~\cite{Sil}, whereas the computations from $\delta_1,\delta_2,\gamma_1,\gamma_2-\gamma_1$ gives a matrix which is different from the one given by ~\cite{Sil}.
\end{remark}

\subsubsection{Example when \texorpdfstring{$\lambda = 2$}{lambda = 2}}\label{sec:lam2}
Fix $\lambda = 2$, in Table~\ref{tab:newlam2} we demonstrate the convergence to the Riemann matrix for levels $0$ through $7$ as defined in Definition~\ref{def:levels}. Indeed, note that the accuracy of the matrix entries increases by about 1 digit with each additional level, resulting in accuracy up to $10^{-5}$ in level 7. The computation was unable to finish on the 8th level. Though the linear equations are very sparse (on the order of 4-5 nonzero coefficients each), due to the size of the system we were not able to push beyond level 8 for the computation. However, we believe there are ways to increase the efficiency of the computation, which may be worth attempting in future work. %We also remark that our experiments coincide with the results following the algebraic computations of~\cite{Sil}. 

\subsubsection{An irrational \texorpdfstring{$\lambda$}{lambda}} \label{sec:irrlam}
In this case, we first approximate the $L$ surface up to a fixed tolerance via continued fractions. As we decrease the tolerance, the denominator of the continued fraction approximation grows, creating finer and finer quadrangulations as the size of the squares is dependent on the size of the denominator (See Table~\ref{tab:lamsqrt}). We selected the value $\lambda = \frac{1+\sqrt{3}}{2}$ since the underlying algebraic curve is defined over the field of rational numbers, given in \cite{Sil}, by $y^2 = x(x^2-1)(x-2)(x-\frac{1}{2}).$

Finally to demonstrate the further convergence, we fix a continued fraction approximation of $\frac{10864}{7953}$ (See Table~\ref{tab:lamsqrt2}). Since the $0$th level already includes squares of size $\frac{1}{2(7953)}$, we were only able to run subdivisions of level 1 and level 2. For comparison we then include the numerical approximations of the Riemann matrices for the continued fraction as well as the original value of $\lambda$. In this case, similarly to the $\lambda = 2$ case, we find that the entries of the matrix are accurate to $10^{-6}$, both for the matrix for the continued fraction and for the original irrational $\lambda$. The entries of the matrices for the continued fraction and for $\lambda$ coincide up to 8 digits.

%\subsubsection{Experiments and Computational Complexity}

%\section{Seeing Convergences of JS-differentials from Igusa's equations}

%\section{$M$-curves}

%\section{How limits work}
%The goal of this section is to use discretization in an effective and explicit way in order to get a limit. 

%\section{Surfaces with many symmetries}
%    The goal here is to look at discrete approximations for $S_6^1 = S_6^2$ and $S_5^2$ as in \cite{DM21}. In particular note that $S_5^2$ is the translation surface formed by identifying opposite sides of the decagon. According to \cite[top of page 6]{McMulleng2} the decagon is generated by a multiple of the decagon form $\omega = \frac{dx}{y}$ on $y^2 = x (x^5 - 1)$. By \cite[Theorem 8.2]{McMullendeca} the decagon is equivalent up to $SL(2,R)$ to the prototypical form of type $(0,1,1,-1)$ with width $t= (2+ \phi)/5$ where $\phi = \frac{1+ \sqrt{5}}{2}$ is the golden ratio.
\subsection{The Jenkins--Strebel representatives} \label{sec:JS}
Given an integer $g \geq 2$, the goal of this section is to use discrete approximations to estimate the curve underlying a Riemann surface of genus $g$ for which we do not a priori know the underlying algebraic curve. Namely we will define $J_g$ to be a square tiled Jenkins--Strebel (JS) representative of the principal stratum. We construct discrete approximations in low genera. Our construction of JS representatives for a holomorphic one-form with one cylinder follows the work of \cite{ZorichJS}.

\subsubsection{Constructing a Jenkins--Strebel representative} 
%= J_g(1,1)
In this section we define the surface $J_g$ for $g
\geq 2$, with some background on how to obtain a Jenkins--Strebel representative of the principal stratum for each genus $g$ and select basis curves for homology. Example~\ref{ex:JSgenus2} illustrates the surface when $g=2$. In general we will construct a Jenkens--Strebel representative by taking a $1\times (4g-4)$ rectangle, identifying sides by translation, where the top and bottom are identified via a permutation $\pi_{J_g}$ of the $4g-4$ horizontal edges and the single vertical side is identified. By varying the side lengths, we can move from $J_g$ through a family of curves which are all Jenkins--Strebel differentials of the same genus. For the sake of tractable computations, we allow ourselves to vary only two parameters that label the vertical edge and the far left horizontal edge in our experiments. We insert the parameters in the parenthesis as $J_g(\lambda,\mu)$ if needed.

We say that $\pi_{J_g}$ is the \emph{permutation associated to} $J_g$.

\begin{theorem}[\cite{ZorichJS} Proposition 2] \label{thm:JSthm}
    Given a genus $g\geq 2$, the Jenkins--Strebel representative associated to the principal stratum composed of unit squares is given by the permutation on $4g-4$ elements where for $k= 1,\ldots, 4g-4$,
    $$\pi_{J_g}(k) = \begin{cases}
        k+1, & k\text{ is odd},\\
        k-1, & k\equiv 2 \pmod{4},\\
        k+3\pmod{4g-4},& k \equiv 0 \pmod 4.
    \end{cases}$$
\end{theorem}
\begin{remark}
    We apply a shear by $\displaystyle{\begin{bmatrix} 1 & -1/2 \\ 0 &1\end{bmatrix}}$, and then perform a cut and paste operation moving the left square over to the right, and then relabel to obtain Theorem~\ref{thm:JSthm} from Proposition 2 of \cite{ZorichJS}
\end{remark}
\begin{remark}\label{rem:4g-4zeros}
    Notice that there are always $4g-4$ elements in the permutation, which comes from the fact that at each singularity, there is an angle of $4\pi$, so there are at most $4$ parallel saddle connections. This gives a total of $4(2g-2)$ possible saddle connections, but since each saddle connection has an incoming and outgoing direction, we have double counted. Thus there are $2(2g-2) = 4g-4$ possible parallel saddle connections. 
\end{remark}

\begin{figure}[htbp]
\centering
\includegraphics[scale=.4]{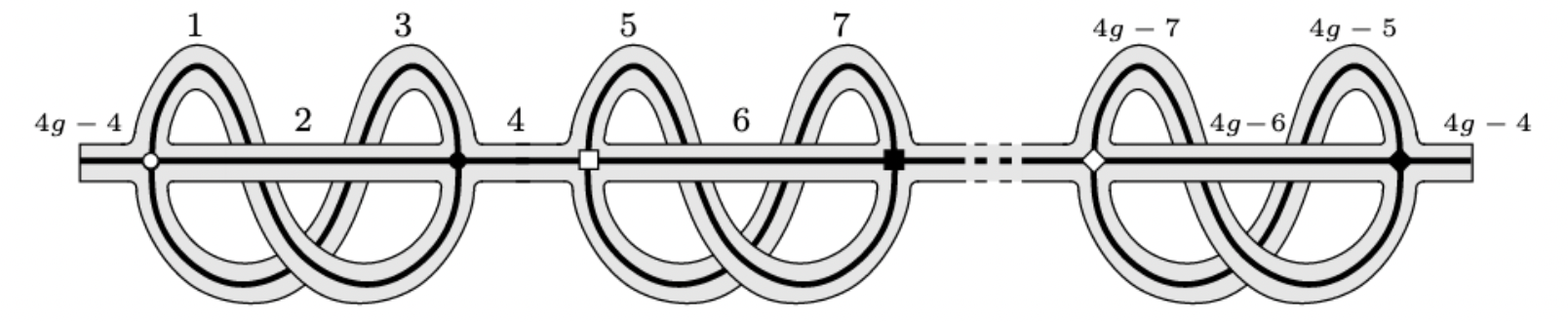}
\caption{Above is the ribbon graph associated to $J_g$. The first pretzel is given by the $1$,$2$, and $3$ edges, and then the two half edges $4g-4$ and $4$. There are a total of $g-1$ pretzels.}
    \label{fig:gribbon}
\end{figure}

Each JS differential is associated to a \textit{ribbon graph}. Indeed we take the graph with vertices given by the $2g-2$ singularities, and edges labelled by the $4g-4$ parallel saddle connections. We then preserve the topological information by contracting the surface with boundary that follows the graph (c.f. Figure~\ref{fig:gribbon}). The fact that we can do this on $J_g$ comes from the fact that JS surfaces have closed horizontal leaves.

For any $g$, we obtain the ribbon graph by gluing $g-1$ \textit{pretzels} together. We define a pretzel to be a ribbon graph on two vertices, with 3 edges connecting the two vertices, and an open edge on each vertex. These pretzels are formed by the fact that $4\pi$ angle allows for a graph of degree 4 at each vertex. Then there are 3 parallel saddle connections between any two vertices, and the $4$th at each pair of saddle connections is used to connect to other pretzels.

The ribbon graph has a top and bottom surrounding each edge. Following around the top edge starting at the top left in Figure~\ref{fig:gribbon}, we start at $4g-4$, then 1, then travel to 2,3,4, and continue in order. Now start at the bottom edge of $4g-4$ on the left. The permutation is now $3, 2, 1,4,7,6,5,$ and so on. For a concrete example, one can perform this exercise with $g=3$ and edge numbers $1-8$ to obtain the permutation as given in Figure~\ref{fig:js3_polygon}.

\begin{figure}[htb]
     \centering
     \begin{subfigure}[b]{\textwidth}
         \begin{center}
             \begin{tikzpicture}[line cap=round,line join=round,x=1.8cm,y=1.8cm]

\draw[-] (0,0) -- (2,0)--(8,0)--(8,1)--(0,1)--(0,0);
\draw (0,.5) node[anchor=east] {$0$};
\draw (8,.5) node[anchor=west] {$0$};

\draw (.5,1) node[anchor=south] {$1$};
\draw (1.5,1) node[anchor=south] {$2$};
\draw (2.5,1) node[anchor=south] {$3$};
\draw (3.5,1) node[anchor=south] {$4$};
\draw (4.5,1) node[anchor=south] {$5$};
\draw (5.5,1) node[anchor=south] {$6$};
\draw (6.5,1) node[anchor=south] {$7$};
\draw (7.5,1) node[anchor=south] {$8$};

\draw (.5,0) node[anchor=north] {$2$};
\draw (1.5,0) node[anchor=north] {$1$};
\draw (2.5,0) node[anchor=north] {$4$};
\draw (3.5,0) node[anchor=north] {$7$};
\draw (4.5,0) node[anchor=north] {$6$};
\draw (5.5,0) node[anchor=north] {$5$};
\draw (6.5,0) node[anchor=north] {$8$};
\draw (7.5,0) node[anchor=north] {$3$};

\fill  [color=blue] (0,0) circle (3pt);
\fill  [color=blue] (1,1) circle (3pt);
\fill  [color=blue] (2,0) circle (3pt);
\fill  [color=blue] (3,1) circle (3pt);
\fill  [color=blue] (8,0) circle (3pt);
\draw  [color=blue] (0,1) circle (3pt);
\draw  [color=blue] (1,0) circle (3pt);
\draw  [color=blue] (2,1) circle (3pt);
\draw  [color=blue] (7,0) circle (3pt);
\draw  [color=blue] (8,1) circle (3pt);
\node[rectangle,draw, color = blue] (r) at (3,0) {};
\node[rectangle,draw, color = blue] (r) at (4,1) {};
\node[rectangle,draw, color = blue] (r) at (5,0) {};
\node[rectangle,draw, color = blue] (r) at (6,1) {};
\node[rectangle,draw, color = blue, fill = blue] (r) at (4,0) {};
\node[rectangle,draw, color = blue, fill = blue] (r) at (5,1) {};
\node[rectangle,draw, color = blue, fill = blue] (r) at (6,0) {};
\node[rectangle,draw, color = blue, fill = blue] (r) at (7,1) {};

\draw[->](.5,0)--(1.5,1);
\draw[color=black] (.5,.25) node {$\alpha_1$};
\draw[->] (4.5,0)--(5.5,1);
\draw[color=black] (4.5,.25) node {$\alpha_2$};
\draw[<-] (6.5,0)--(7.5,1);
\draw[color=black] (6.5,.25) node {$\alpha_3$};
\draw[<-] (2.5,0)--(3.5,1);
\draw[color=black] (2.5,.25) node {$\alpha_3$};
\draw[->](1.5,0)--(.5,1);
\draw[color=black] (1.5,.25) node {$\beta_1$};
\draw[->] (5.5,0)--(4.5,1);
\draw[color=black] (5.5,.25) node {$\beta_2$};
\draw[->, dashed] (2.5,1)--(3.5,0);
\draw[color=black] (7.5,.25) node {$\beta_3$};
\draw[->, dashed](6.5,1) --(7.5,0);
\draw[color=black] (3.5,.25) node {$\beta_3$};

%\draw[color=black] (-0.39,-0.22) node {$n$};

\end{tikzpicture}
         \end{center}
         \caption{Here is the polygonal representation of the JS surface for $g=3$. The associated permutation is $\pi_{J_3} = \begin{pmatrix}1&2&3&4&5&6& 7 & 8\\
    2&1&4&7&6&5&8&3\end{pmatrix}$. The sides which are not $3\pmod{4}$ are labelled by the homology basis curve that is used for the periodicity equations in Algorithm~\ref{alg:general}. The curve $\beta_3$ is represented by a dashed line since more care must be taken to construct the perioidicity equations. (See Algorithm~\ref{sec:AlgForJS}.)}
         \label{fig:js3_polygon}
     \end{subfigure}
     \begin{subfigure}[b]{\textwidth}
         \centering
        \begin{tikzpicture}[scale = 1.5]
\draw[smooth] (0,1) to[out=30,in=150] (2,1) to[out=-30,in=210] (3,1) to[out=30,in=150] (5,1) to[out=-30,in=210] (6,1) to[out=30,in=150] (8,1) to[out=-30,in=30] (8,-1) to[out=210,in=-30] (6,-1) to[out=150,in=30] (5,-1) to[out=210,in=-30] (3,-1) to[out=150,in=30] (2,-1) to[out=210,in=-30] (0,-1) to[out=150,in=-150] (0,1);
\draw[smooth] (0.4,0.1) .. controls (0.8,-0.25) and (1.2,-0.25) .. (1.6,0.1);
\draw[smooth] (0.5,0) .. controls (0.8,0.2) and (1.2,0.2) .. (1.5,0);
\draw[smooth] (3.4,0.1) .. controls (3.8,-0.25) and (4.2,-0.25) .. (4.6,0.1);
\draw[smooth] (3.5,0) .. controls (3.8,0.2) and (4.2,0.2) .. (4.5,0);
\draw[smooth] (6.4,0.1) .. controls (6.8,-0.25) and (7.2,-0.25) .. (7.6,0.1);
\draw[smooth] (6.5,0) .. controls (6.8,0.2) and (7.2,0.2) .. (7.5,0);
\draw[color=red] (4.0,0.15) arc(270:90:0.3 and 1.14/2) node[sloped,pos=.5,allow upside down]{\arrowIn}node[above, pos=0] {$\alpha_3$};;
\draw[dashed, color=red] (4.0,0.15) arc(270:450:0.3 and 1.14/2);
\draw[color=red] (1.0,0.15) arc(270:90:0.3 and 1.14/2) node[sloped,pos=.5,allow upside down]{\arrowIn}node[above, pos=0] {$\alpha_1$};;
\draw[dashed, color=red] (1.0,0.15) arc(270:450:0.3 and 1.14/2);
\draw[color=red] (7.0,0.15) arc(270:90:0.3 and 1.14/2) node[sloped,pos=.5,allow upside down]{\arrowIn}node[above, pos=0] {$\alpha_2$};;
\draw[dashed, color=red] (7.0,0.15) arc(270:450:0.3 and 1.14/2);
\draw[color=cyan] (2.0,0.0) arc(0:360:1 and 1.14/2) node[sloped,pos=.5,allow upside down]{\arrowIn}node[above, pos=.75] {$\beta_1$};;
\draw[color=cyan] (5.0,0.0) arc(0:360:1 and 1.14/2) node[sloped,pos=.5,allow upside down]{\arrowIn}node[above, pos=.75] {$\beta_3$};;
\draw[color=cyan] (8.0,0.0) arc(0:360:1 and 1.14/2) node[sloped,pos=.5,allow upside down]{\arrowIn} node[above, pos=.75] {$\beta_2$};;

\draw[thick] (5.5,0.5) arc(0:180:1.5 and 1.14/3) node[sloped,pos=.5,allow upside down]{\arrowIn}node[above, pos=.4] {$8$};;
\draw[thick] (2.5,-0.5) arc(180:360:1.5 and 1.14/3) node[sloped,pos=.5,allow upside down]{\arrowIn}node[below, pos=.4] {$4$};;
\draw[smooth, thick] (2.5,0.5) .. controls (-1,2) and (-1,-2) .. (2.5,-.5)node[sloped,pos=.3,allow upside down]{\arrowIn}node[above, pos=.3] {$2$};;
\draw[smooth, thick] (5.5,-.5) .. controls (9.5,-2) and (9,2) .. (5.5,.5)node[sloped,pos=.7,allow upside down]{\arrowIn}node[above, pos=.71] {$6$};;
%%%%1
\draw[smooth, thick] (1.5,0) .. controls (1.5,-.25) and (2,-.5) .. (2.5,-.5)node[sloped,pos=.3,allow upside down]{\arrowIn}node[above, pos=.3] {$1$};;
\draw[smooth, thick] (2.5,.5) .. controls (2.4,.75) and (2.3,.8) .. (2.2,.9) node[sloped,pos=.5,allow upside down]{\arrowIn}node[above, pos=1] {$1$};;
\draw[smooth, dashed] (2.15,.9) .. controls (1.9,.9) and (1.6,.5) .. (1.5,0);

%%% 3
\draw[smooth, thick] (2.5,.5) .. controls (2.75,.5) and (3,.5) .. (3.5,0)node[sloped,pos=.5,allow upside down]{\arrowIn}node[above, pos=.71] {$3$};;
\draw[smooth, thick] (2.7,-.875) .. controls (2.65,-.8) and (2.45,-.6) ..(2.5,-.5)  node[sloped,pos=.25,allow upside down]{\arrowIn}node[below, pos=.1] {$3$};;
\draw[smooth, dashed] (3.5,0) .. controls (3.5,-.5) and (3,-.75) .. (2.75,-.875);

%%% 7
\draw[smooth, thick] (5.5,-.5) .. controls (5,-.5) and (4.5,-.25)  .. (4.5,0) node[sloped,pos=.25,allow upside down]{\arrowIn}node[above, pos=.2] {$7$};;
\draw[smooth, thick] (5.2,.9) .. controls (5.3,.8) and (5.4,.75) ..(5.5,.6)  node[sloped,pos=.5,allow upside down]{\arrowIn}node[above, pos=0] {$7$};;
\draw[smooth, dashed] (5.2,.9) .. controls (4.9,.9) and (4.6,.5) .. (4.5,0);

%%% 5
\draw[smooth, thick] (6.5,0) .. controls (6,.5) and (5.75,.5) .. (5.5,.5)node[sloped,pos=.5,allow upside down]{\arrowIn}node[above, pos=.3] {$5$};;
\draw[smooth, thick] (5.5,-.5) .. controls (5.65,-.8) and (5.45,-.6) .. (5.7,-.875) node[sloped,pos=.9,allow upside down]{\arrowIn}node[below, pos=1] {$5$};;
\draw[smooth, dashed] (6.5,0) .. controls (6.5,-.5) and (6,-.75) .. (5.75,-.875);

%%% 0
\draw[smooth, thick] (2.5,.5) .. controls (2.6,.75) and (2.5,.8) .. (2.5,.85) node[sloped,pos=.5,allow upside down]{\arrowIn}node[above, pos=1] {$0$};;
\draw[smooth, thick] (2.25,-.875) .. controls (2.45,-.8) and (2.5,-.6) ..(2.5,-.5)  node[sloped,pos=.25,allow upside down]{\arrowIn}node[below, pos=.1] {$0$};;
\draw[smooth, dashed] (2.5,.85) .. controls (2,.5) and (2,-.5) .. (2.25,-.875);

\filldraw[color=blue, fill=white] (2.5,.5) circle(2pt);
\filldraw[color=blue] (2.5,-.5) circle(2pt);
\filldraw[color=blue, fill=white] (5.5,-.45) rectangle (5.65,-.6);
\filldraw[color=blue] (5.5,.45) rectangle (5.65,.6);
\end{tikzpicture}
         \caption{Here is the topological representation of the polygonal surface for $g=3$ keeping track of edge identifications and zeroes of the $1$-form, and the intersections of the edges with the standard homology basis.}
         \label{fig:js3_topological}
     \end{subfigure}
     
        \caption{The polygonal representation and topological representation giving information of how homology vectors behave under identified sides of the JS surface for genus 3.}
        \label{fig:JS3}
\end{figure}

The next step we must take is to fix a basis of homology to help in determining the periodicity equations. We first consider the topological picture, as in Figure~\ref{fig:js3_topological}. We fix a standard symplectic homology basis curves $\alpha_j, \beta_j$ for $j=1,\ldots, g$ where the pairs $\alpha_i,\beta_i$ are oriented to have intersection number 1, and all other pairs have intersection zero. For higher genus, the picture is best seen when $\alpha_g,\beta_g$ are in the center, and all other attached tori are equally spread out and only attached to the $\alpha_g,\beta_g$ torus. The idea is each pretzel in the ribbon graph adds one torus and two zeros. Notice that the symplectic basis is less obvious in the polygonal picture. In Figure~\ref{fig:js3_polygon}, by removing the second reference to $\alpha_3$ used for periodicity identifications, the set of curves are consistent in the intersection numbers to be a symplectic basis, and one can see the curves are not homologous by cutting the  surface into $g$ tori with boundary components formed by cutting vertically along the zeroes. The curves $\alpha_j, \beta_j$ then form the standard two symplectic basis vectors of homology of the torus excluding the basis curves relative to the boundary. 

Next mark the $2g-2$ zeroes of the one-form, between the $\alpha_g,\beta_g$ torus and the corresponding numbers of that pretzel.  For example in Figure~\ref{fig:js3_topological}, the filled and unfilled circles belong to the first pretzel, and the filled and unfilled squares belong to the second pretzel. Next label each edge coming out of a vertex in the correct order. To determine the correct order, for example in Figure~\ref{fig:js3_polygon}, consider the unfilled circle. Starting at the $0$ side, we travel in a circle with the unfilled circle vertex to our left. From side $0$, we go to side $1$, keep traveling with the circle on our left, cross side $2$, then $3$, and $8$ before returning to side zero. This order matches the order seen in Figure~\ref{fig:js3_topological}.

Now to connect edges, we have 3 cases. First we select the $1,2\pmod{4}$ edges in the $j$th pretzel to be the edges crossing $\beta_j$, $\alpha_j$, respectively. Next the $3\pmod{4}$ edges cross the $\beta_g$. Third, the $0 \pmod{4}$ sides are either $0$ and cross no homology curves, or cross $\alpha_g$. In this manner, every homology curve is crossed by at least $1$ edge. We then mark the respective homology curves and their directions in the polygonal picture by keeping track of directions. For example in \ref{fig:js3_topological}, as we travel along side $2$ from open circle to filled circle $\alpha_1$ crosses from left to right. Similarly $\alpha_3$ crosses from left to right as we travel along side $8$ from filled square to open circle, but this is exactly the opposite direction as side $2$ in the polygonal picture Figure~\ref{fig:js3_polygon}, giving the marked directions of $\alpha_1$ and $\alpha_3$.

\begin{remark}
    Notice the choice of the principal stratum is for simplicity. To construct other orders of zeroes, \cite{ZorichJS} combines zeroes in the principal stratum. When the zeroes are combined, the same method of carefully following edge identifications in the polygon and how they connect to the homology basis works as well.
\end{remark}

\begin{figure}[htbp]
    \centering
   \begin{tikzpicture}[line cap=round,line join=round,x=1 cm,y=1 cm, scale = .9]
\draw[-] (0,0) --(16,0)--(16,1)--(0,1)--(0,0);

\draw[-] (2,2) -- (14,2)--(14,3)--(2,3)--(2,2);
\begin{scriptsize}

\draw (0,.5) node[anchor=east] {$0$};
\draw (16,.5) node[anchor=west] {$0$};
\draw (2,2.5) node[anchor=east] {$0$};
\draw (14,2.5) node[anchor=west] {$0$};

\draw (.5,1) node[anchor=south] {$1$};
\draw (1.5,1) node[anchor=south] {$2$};
\draw (2.5,1) node[anchor=south] {$3$};
\draw (3.5,1) node[anchor=south] {$4$};
\draw (4.5,1) node[anchor=south] {$5$};
\draw (5.5,1) node[anchor=south] {$6$};
\draw (6.5,1) node[anchor=south] {$7$};
\draw (7.5,1) node[anchor=south] {$8$};
\draw (8.5,1) node[anchor=south] {$9$};
\draw (9.5,1) node[anchor=south] {$10$};
\draw (10.5,1) node[anchor=south] {$11$};
\draw (11.5,1) node[anchor=south] {$12$};
\draw (12.5,1) node[anchor=south] {$13$};
\draw (13.5,1) node[anchor=south] {$14$};
\draw (14.5,1) node[anchor=south] {$15$};
\draw (15.5,1) node[anchor=south] {$16$};

\draw (.5,0) node[anchor=north] {$2$};
\draw (1.5,0) node[anchor=north] {$1$};
\draw (2.5,0) node[anchor=north] {$4$};
\draw (3.5,0) node[anchor=north] {$7$};
\draw (4.5,0) node[anchor=north] {$6$};
\draw (5.5,0) node[anchor=north] {$5$};
\draw (6.5,0) node[anchor=north] {$8$};
\draw (7.5,0) node[anchor=north] {$11$};
\draw (8.5,0) node[anchor=north] {$10$};
\draw (9.5,0) node[anchor=north] {$9$};
\draw (10.5,0) node[anchor=north] {$12$};
\draw (11.5,0) node[anchor=north] {$15$};
\draw (12.5,0) node[anchor=north] {$14$};
\draw (13.5,0) node[anchor=north] {$13$};
\draw (14.5,0) node[anchor=north] {$16$};
\draw (15.5,0) node[anchor=north] {$3$};

\draw[->](.5,0)--(1.5,1);
\draw[color=black] (.5,.25) node {$\alpha_1$};
\draw[->] (4.5,0)--(5.5,1);
\draw[color=black] (4.5,.25) node {$\alpha_2$};
\draw[->] (8.5,0)--(9.5,1);
\draw[color=black] (8.5,.25) node {$\alpha_3$};
\draw[->] (12.5,0)--(13.5,1);
\draw[color=black] (12.5,.25) node {$\alpha_4$};
\draw[<-] (14.5,0)--(15.5,1);
\draw[color=black] (14.5,.25) node {$\alpha_5$};
\draw[<-] (10.5,0)--(11.5,1);
\draw[color=black] (10.5,.25) node {$\alpha_5$};
\draw[<-] (6.5,0)--(7.5,1);
\draw[color=black] (6.5,.25) node {$\alpha_5$};
\draw[<-] (2.5,0)--(3.5,1);
\draw[color=black] (2.5,.25) node {$\alpha_5$};
\draw[->](1.5,0)--(.5,1);
\draw[color=black] (1.5,.25) node {$\beta_1$};
\draw[->] (5.5,0)--(4.5,1);
\draw[color=black] (5.5,.25) node {$\beta_2$};
\draw[->] (9.5,0)--(8.5,1);
\draw[color=black] (9.5,.25) node {$\beta_3$};
\draw[->] (13.5,0)--(12.5,1);
\draw[color=black] (13.5,.25) node {$\beta_4$};
\draw[->, dashed] (2.5,1)--(3.5,0);
\draw[->, dashed] (6.5,1)--(7.5,0);
\draw[->, dashed] (10.5,1)--(11.5,0);
\draw[->, dashed] (14.5,1)--(15.5,0);

%\draw[color=black] (-0.39,-0.22) node {$n$};
%%%%%%%%%%%%%%%%%%%%%%%%%%%%%%%%%%%%%%%%

\draw (2.5,3) node[anchor=south] {$1$};
\draw (3.5,3) node[anchor=south] {$2$};
\draw (4.5,3) node[anchor=south] {$3$};
\draw (5.5,3) node[anchor=south] {$4$};
\draw (6.5,3) node[anchor=south] {$5$};
\draw (7.5,3) node[anchor=south] {$6$};
\draw (8.5,3) node[anchor=south] {$7$};
\draw (9.5,3) node[anchor=south] {$8$};
\draw (10.5,3) node[anchor=south] {$9$};
\draw (11.5,3) node[anchor=south] {$10$};
\draw (12.5,3) node[anchor=south] {$11$};
\draw (13.5,3) node[anchor=south] {$12$};

\draw (2.5,2) node[anchor=north] {$2$};
\draw (3.5,2) node[anchor=north] {$1$};
\draw (4.5,2) node[anchor=north] {$4$};
\draw (5.5,2) node[anchor=north] {$7$};
\draw (6.5,2) node[anchor=north] {$6$};
\draw (7.5,2) node[anchor=north] {$5$};
\draw (8.5,2) node[anchor=north] {$8$};
\draw (9.5,2) node[anchor=north] {$11$};
\draw (10.5,2) node[anchor=north] {$10$};
\draw (11.5,2) node[anchor=north] {$9$};
\draw (12.5,2) node[anchor=north] {$12$};
\draw (13.5,2) node[anchor=north] {$3$};

\draw[->](2.5,2)--(3.5,3);
\draw[color=black] (2.5,2.25) node {$\alpha_1$};
\draw[->] (6.5,2)--(7.5,3);
\draw[color=black] (6.5,2.25) node {$\alpha_2$};
\draw[->] (10.5,2)--(11.5,3);
\draw[color=black] (10.5,2.25) node {$\alpha_3$};
\draw[<-] (12.5,2)--(13.5,3);
\draw[color=black] (12.5,2.25) node {$\alpha_4$};
\draw[<-] (8.5,2)--(9.5,3);
\draw[color=black] (8.5,2.25) node {$\alpha_4$};
\draw[<-] (4.5,2)--(5.5,3);
\draw[color=black] (4.5,2.25) node {$\alpha_4$};
\draw[->](3.5,2)--(2.5,3);
\draw[color=black] (3.5,2.25) node {$\beta_1$};
\draw[->] (7.5,2)--(6.5,3);
\draw[color=black] (7.5,2.25) node {$\beta_2$};
\draw[->] (11.5,2)--(10.5,3);
\draw[color=black] (11.5,2.25) node {$\beta_3$};
\draw[->, dashed] (8.5,3)--(9.5,2);
\draw[->, dashed] (12.5,3)--(13.5,2);
\draw[->, dashed] (4.5,3)--(5.5,2);
%\draw[color=black] (13.5,2.25) node {$\beta_4$};

\end{scriptsize}
\end{tikzpicture}
    \caption{Surface representations of $J_4$ (above) and $J_5$ (below) with marked homology basis curves, the basis curve $\beta_g$  is represented by the dashed line.}
    \label{fig:g45JS}
\end{figure}
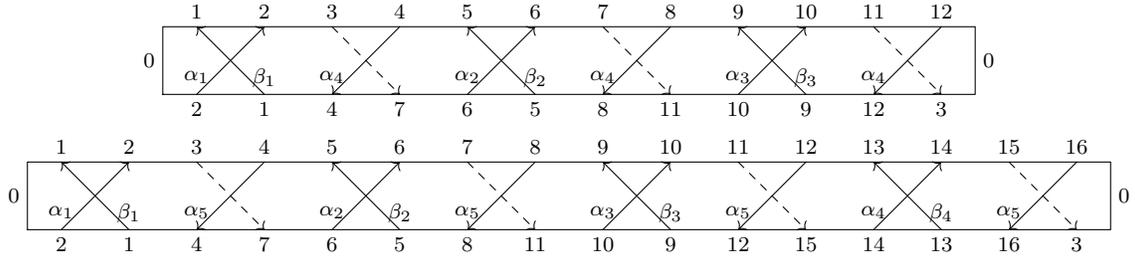

In addition to computations for the two examples given in Figure~\ref{fig:g=2JS} and Figure~\ref{fig:JS3}, we will also give examples in genus $4$ and $5$ (Figure~\ref{fig:g45JS}), where the choice of homology basis follows the same strategy outlined above.

\subsubsection{Algorithm for discrete approximation of JS surfaces}\label{sec:AlgForJS}
Let $g$ be the genus of the JS representative. Let $n \in \mathbb{N}\cup{\{0\}}$ be the level of approximation. For simplicity of the description we here fix all side lengths to be 1, but implementation with varying side lengths can be found on our MathRepo page~\cite{CelFaiMan}. We will note here, however, that the algorithm requires side lengths to be rational numbers, with odd numerator and denominator, due to the bicolored structure. 

\begin{enumerate}
    \item \textbf{Constructing an initial bipartite quadrangulation.} Since all the side lengths are 1, the shape is already divided into squares. In order to bi-color the square tiling and maintain the vertex identifications, there must be an even number of squares on each side length, which is always true since $4g-4$ is always even. Hence we define the step size to be $1$.
       
    \item \textbf{Quadrangulations for further levels of approximation.} Each square of side length $1$ will be divided into $3^{2n}$ squares, and so the $n$th level approximation will consist of squares of size ${1}/{3^n}.$ 
    \item \textbf{Labelling vertices.} We label the vertices with the following bounds:
    $$x_{i,j}= \left(\frac{i}{3^n}, \frac{j}{3^n } \right) \quad \text{ for } \,0 \leq i \leq 3^n(4g-4)\, \text{ and }\, 0\leq j\leq 3^n.$$

    \item \textbf{Holomorphicity equations.} For each bottom left of a square, we have a new holomorphicity equation. So in this case, we have a total of $3^{2n} (4g   - 4)$ equations with indices given by $ 0\leq i \leq  3^n(4g-4)-1$ and $0\leq j\leq 3^n -1$.
  
    \item \textbf{Periodicity equations.} The basis of homology was chosen so that each edge aside from the 0 edge and the edges which are congruent to $3 \pmod{4}$  is identifed via a single basis of homology vector. For all of the following periods, the parity $p$ is determined by $p = b$ if $i+j \equiv 0 \pmod 2$ and $p=w$ otherwise, so we only need to determine the periodicity relationships.
    \begin{itemize}
        \item \textbf{$\mathbf{A_{k+1}}$ periods for $\mathbf{k=0,\ldots, g-2}$.} We have $(g-1)(3^n+1)$ equations where for each $k$ there are $3^n+1$ equations, and the basis vectors occur at $2\pmod{4}$ to give
        $$A_{k+1}^p = x_{i+3^{n} , 3^n} - x_{i,0}\quad \text{for } \,\,  4 k 3^n \leq i\leq 3^n (4k+1) .$$
        \item \textbf{$\mathbf{A_g}$ periods.} We have $(g-1)(3^n+1)$ equations where there are $3^n+1$ equations occuring at each of the $0 \pmod{4}$ sides given by
        $$A_{g}^p = x_{i,0} - x_{i+3^{n} , 3^n}\quad \text{for } \,\, 3^n(4k+2)\leq i\leq 3^n(4k+3). $$
        \item \textbf{$\mathbf{B_k}$ periods for $\mathbf{k= 1,\ldots, g-1}$.} We use the same indexing as the $A_{k+1}$ periods to have $(g-1)(3^n+1)$ equations given by
        $$B_{k+1}^p = x_{i, 3^n} - x_{i+3^{n} , 0} \quad \text{for } \,\,  4 k 3^n \leq i\leq 3^n (4k+1). $$
        \item \textbf{$\mathbf{B_g}$ periods and $\mathbf {0 \pmod{4}}$ side identifications.} The homology curve is crossed by all of the sides $3\pmod{4}$ , and the side $0$ crosses no homology curve. To find the correct identifications, we construct relations from highest numbers to lowest numbers. The idea is for each unidentified edge, we follow all paths around vertices which identify the two possible sides, except for the side $3$ which we choose to be the curve we follow around $\beta_3$. We go through this carefully following the identification in the example of $g=3$ following the image in Figure~\ref{fig:js3_topological}. We refer to the top and bottom of every edge by looking at the top and bottom in the polygonal representation of Figure~\ref{fig:js3_polygon}. 
        \begin{itemize}
            \item \textbf{Sides $\mathbf{3\pmod{4}}$ bigger than $\mathbf{3}$.} On side $7$, we describe the path travelling from $7_{bottom}$ to $7_{top}$, by travelling around the filled square to the right. To do this we cross side $6$ from bottom to top, which is identified via $\alpha_2$. Next we cross side $5$ top to bottom giving $-\beta_2$, and finally side $8$ bottom to top giving a $-\alpha_3$. So all together,
            $$7_{top}-7_{bottom} = \alpha_2 - \beta_2 - \alpha_3.$$
            Following along a circle with the unfilled square to the left verifies the choice that side $4$ is identified by $\alpha_3.$ In general working with $3\pmod{4}$ sides not equal to 3 we have $(g-2) (3^n+1)$ equations of the form
        $$x_{(i+3\cdot3^n), 3^n} -x_{i , 0} = A_{k+2} - B_{k+2} -A_g  \quad \text{for }  \begin{cases} 0 \leq k \leq g-3,\\
        (4k+3)3^n \leq i\leq 3^n (4k+4).\end{cases}$$
            \item \textbf{Side $\mathbf{3}$.} On the $3$ side, travelling around the zeroes as in the $7$ sides will be used to determine the identifications of the $0$ sides. We choose the $3$ side to encode the crossing information of $\beta_g$ since it is included in all genus $g\geq 2$. When $g=2$, the curve $\beta_2$ only is crossed by side $3$. However for $g=3$, we may follow the image in Figure~\ref{fig:js3_topological}. Starting at the bottom of $3$, we travel in the direction of $-\beta_3$, crossing side $7$ from top to bottom before we complete $\beta_3$ to return to the bottom of $3$. Thus 
            $$3_{top}-3_{bottom} = -\beta_3 + 7_{bottom}-7_{top} = -\beta_3 + \alpha_3 - \alpha_2 + \beta_2.$$
            For higher genus, there are $g-2$ total $3\pmod{4}$ sides each contributing $\alpha_3$ and some $\beta_k - \alpha_k$.
            
            Setting $$\iota = i+3 \cdot 3^n \pmod{3^n(4g-4)}$$ we have $3^n+1$ equations for identifying the top and bottom of side $3$ given by 
        $$x_{\iota , 3^n} -x_{i , 0} = -B_g + (g-2)A_g + \sum_{j=2}^{g-1}(B_j-A_j)  \quad \text{for }  (4g-5)3^n \leq i\leq 3^n (4g-4) .$$
            \item \textbf{Side $\mathbf{0}$.} We conclude by identifying the zero sides. We will travel around with the open circle on the left, and leave it to the reader to verify the same result holds for travelling around the closed circle to the right. We start at the left side of zero, crossing sides $1$, $2$, $3$, $8$. This gives
            $$0_{right}-0_{left} = -\beta_1  + \alpha_1 + 3_{bottom}-3_{top} - \alpha_3 =  \beta_3 -2\alpha_3 + \alpha_2 - \beta_2 + \alpha_1-\beta_1.$$
            In general, we always follow the same crossing pattern, of $-\beta_1 + \alpha_1 - \alpha_3$, and then this must be combined with the information about the $3$ side. This gives $3^n+1$ equations of the form
            $$x_{3^n(4g-4),j}- x_{0,j} = B_g - (g-1)A_g +\sum_{j=1}^{g-1} (B_j-A_j) \quad \text{for } \,\, 0\leq j\leq 3^n.$$
            
        \end{itemize}

        \begin{comment}We first have $3^n(g-2)$ identifications when the side is $0\pmod{4}$, where $\beta_g$ goes from bottom to top, giving equations for $k=0,\ldots, g-3$
        $$B_{g}^p = x_{i+ 3^n , 3^n} - x_{i, 0} \quad \text{for }  (4 k+2) 3^n \leq i\leq 3^n (4k+3) .$$
        In addition, there are $3^n(g-1)$ identifications for the sides that are $3\pmod{4}$, where $\beta_g$ goes from top to bottom, giving equations for $k=0,\ldots, g-2$
        $$B_{g}^p = x_{i , 0} - x_{i+(4k+3 \pmod{4g-4}) , 3^n} \quad \text{for }  (4 k+3) 3^n \leq i\leq 3^n (4k+4) .$$ \end{comment}
        %\item \textbf{Identifying the 0 sides.} The last components that are not identified?
    \end{itemize}
    \item\textbf{Final Normalizations.} In the final two normalizations we have the following number of equations:
    \begin{itemize}
        \item 2 equations for normalization of holomorphic function.
        \item For each $k=1,\ldots, g$, there are $2g$ equations normalizing to the canonical basis. 
    \end{itemize}
    
    \item \textbf{Solving a system of equations for the discrete approximation.} 
    For the $k$th row of the period matrix with $k=1,\ldots,g$, we obtain the equations by solving the system with:
    \begin{itemize}
        \item Total of $3^{2n}(4g-4) + 3^n(4g-3) + 4g+1$ variables. With $(3^{n}(4g-4)(3^n+1) + (3^n+1)$ variables $x_{i,j}$, and $4g$ variables coming from $A_k^p, B_k^p$ with $k=1,\ldots, g$ and parity given by $b$ and $w$. 
        \item Total of $3^{2n}(4g-4) + 3^n(4g-3) +6g-1$ equations. With $3^{2n}(4g-4)$ holomorphicity equations, $(2g-2)(3^n+1)$ periodicity equations for the $A$ periods, $(2g-2)(3^n+1)$ equations for the $B$ periods and $3\pmod{4}$ sides, $3^n+1$ equations for the $0$ sides, and $2g+2$ normalizing equations.
    \end{itemize}
    Since $g\geq 2$, we have a system of equations over determined by $2g-2$ equations, and these are not conflicting with a unique solution, since they all describe the unique discrete holomorphic differential guaranteed to exist by \cite[Theorem~6.8]{BobGun}.
\end{enumerate}

\subsubsection{Experiments in low genus}\label{sec:g2JSexp}

We consider the JS respresentative for $g=2,3,4,5$ in the principal stratum. We only consider 2-parameter family $J_g(\lambda,\mu)$ as noted in the beginning of Section~\ref{sec:JS}. When $g=2$,  such a surface has been exhibited in Example~\ref{ex:JSgenus2}. Setting the side lengths to be $\lambda= \mu =1$, we present our experiments of the discrete period matrix approximations up to level 7 in Table~\ref{tab:JS2}--\ref{tab:JS5}, which use our algorithm that has been described in Section~\ref{sec:AlgForJS}. Our computations with the approximations in theta functions encourage us to make the following conclusions for low genus.

\noindent {\textbf{Genus $\mathbf{2}$.}} We now take the discrete Riemann matrix of level 7 from Table~\ref{tab:JS2}. As in Example~\ref{ex:SilholLambda2}, we compute the odd theta constants and then approximate the six branch points of the hyperelliptic curve corresponding the $J_2$ surface via the SageMath package~\cite{BruGan}:
 \begin{align*}
      &\alpha_1 := -3.55001177927944 + \mathbf{i}\cdot 9.27369555271397,\\
 &\alpha_2 := -0.0360027110167584 - \mathbf{i}\cdot0.0940498797751955,\\
&\alpha_3 := 0.603906137193071 + \mathbf{i}\cdot3.24517640725254,\\
&\alpha_4 := 0.0554252204362169 - \mathbf{i}\cdot0.297835386410189,\\
&\alpha_5 := 3.90800485599692 - \mathbf{i}\cdot7.79154768793860,\\
&\alpha_6 := 0.0514341663705383 + \mathbf{i}\cdot0.102546382318450.
 \end{align*}

\noindent We observe that these values are pairwise reciprocal, more precisely $\alpha_1\cdot \alpha_2=\alpha_3\cdot \alpha_4=\alpha_5\cdot \alpha_6=1$ up to a numerical round off, which can be sharpened by working with a higher precision complex field. Computing also the 10 even theta constants, we see that some of these values coincide up pairwise to numerical error. These pairs of the constants are in the following pairs of characteristics: 
$$\left\{\begin{bmatrix} 1 & 0 \\ 0 & 0 \end{bmatrix},\begin{bmatrix} 0 & 1 \\ 0 & 0 \end{bmatrix}\right\}, \quad \left\{\begin{bmatrix} 0 & 0 \\ 1 & 0 \end{bmatrix},\begin{bmatrix} 0 & 0 \\ 0 & 1 \end{bmatrix}\right\}, \quad \left\{\begin{bmatrix} 1 & 0 \\ 0 & 1 \end{bmatrix},\begin{bmatrix} 0 & 1 \\ 1 & 0 \end{bmatrix}\right\}.$$
Concluding from similar computations as above, we observe that the reciprocity is respected for some other rational $\lambda$ parameters and $\mu=1$. On the other hand, this phenomenon does not occur if we change $\mu$.   

\noindent Therefore, our experiments suggest making the following conjecture:

\begin{conjecture}
The family of hyperelliptic curves corresponding to the family of the translation surfaces, $J_2(\lambda, 1)$ of genus 2, in the stratum $\mathcal{H}(1,1)$, is given by the equation: 
\begin{equation}
\label{eq:conjCurve}
   y^2=(x-a)(x-1/a)(x-b)(x-1/b)(x-c)(x-1/c) 
\end{equation}
for some complex parameters $a,b,c$. 
\end{conjecture}

 A first remark about the conjecture is that the curve~\eqref{eq:conjCurve} can be transformed to the curve by a projective transformation of $\mathbb{P}^1$: 
%namely mapping $a,b,c$ to $0,1,\infty$,
\begin{equation*}
   y^2=x(x-1)(x-A)(x-B)(x-B/(1-A-B)) 
\end{equation*}
for some complex parameters $A,B$. The model~\eqref{eq:conjCurve} also manifests that this hyperelliptic curve has an extra involution, namely $(x,y)\mapsto(1/x,y)$. This suggests studying of the translation surface $J_2(\lambda,1)$, in particular to reconstruct the underlying algebraic curve via exact computations, akin to the work of \cite{Sil,Rod} for the $L$-shape. We also remark that there are a total of 5 possible choices of side lengths when constructing $J_2$, and we do not expect the symmetry of Weierstrass points to hold for every choice. Indeed by restricting to only allowing $\lambda$ to change, we observe Conjecture~\ref{eq:conjCurve} still holds in the cases tested, whereas allowing $\mu$ to change no longer preserved the symmetries of the fixed points. Indeed the change in symmetry can already be seen on the diagonal entries as shown in Table~\ref{tab:JS2New}, which by continuity, we expect similar results for irrational $\lambda$.

In \cite{rubinstein}, the author is able to construct algebraic equations for genus two curves by approximating the coefficients numerically to high precision with transcendental methods, and then deducing the algebraic numbers using their continued fraction expansions. Once one has the algebraic numbers, it is significantly easier to prove that the equations are correct. In our case, it appears that we need more digits to be able to draw conclusions from the continued fraction expansions. In particular, the methods in section 5 of \cite{rubinstein} do not lead to any clear guesses as to what the algebraic numbers may be. An improvement to our algorithm's implementation, allowing one to compute further levels of subdivision, may lead to these methods becoming applicable.

\noindent {\textbf{Genus $\mathbf{3}$.}} We take the level 7 approximation in Table~\ref{tab:JS3}. Among the 36 even theta constants that we compute in SageMath, one of them gets closer to zero as the precision is increased. As the hyperelliptic curves of genus 3 are characterized with the condition of at least one vanishing even theta constant, we state Conjecture~\ref{con:g3g4} for $g=3$.

\noindent {\textbf{Genus $\mathbf{4}$.}} We first look for evidence that the discrete Riemann matrix estimates a Riemann matrix of an algebraic curve. So we evaluate the discrete Riemann matrix of level 7 from Table~\ref{tab:JS4} in the Schottky-Igusa modular form~\cite{Igusa} as the underlying precision increases. We observe the values approximate to zero. For references see Section~\ref{sec:RS}. Similar to the case above, our computations in SageMath support that there are 10 vanishing even theta constants, which inspires us to state Conjecture~\ref{con:g3g4} for $g=4$.

\noindent {\textbf{Genus $\mathbf{5}$.}} Plugging the discrete period matrix of level 7 from Table~\ref{tab:JS4} into the three equations~\cite[Proposition 1.2]{FarGruMan}, we estimate each of the three values at zero. In addition, we observe that the number of even theta constants that converge the zero is more than 10, which concludes Conjecture~\ref{con:g3g4} for $g=5$. 

\begin{conjecture}\label{con:g3g4}

The surface $J_g(\lambda,1)$ is hyperelliptic for $g=3,4,5$.

\end{conjecture}

One of the primary results of \cite{ZorichJS} was to give JS representatives in every connected component of every stratum, of which some of these connected components are called hyperelliptic. As stated in remark 3 of \cite{KZ03} we note that a hyperelliptic Riemann surface is not always contained in a hyperelliptic connected component. Thus we can expect these specific elements of the principal stratum to be hyperelliptic, but shouldn't expect that we keep a hyperelliptic involution once we allow less symmetries by allowing the $4g-3$ side lengths to be changed.

\bibliographystyle{alpha}
%\nocite{*}
\bibliography{Sources}
\noindent T\"urk\"u \"Ozl\"um \c{C}elik,
\hfill  {\tt turkuozlum@gmail.com}
\vspace{-4pt}

\noindent Samantha Fairchild,
\hfill  {\tt fairchil@mis.mpg.de}
\vspace{-4pt}

\noindent Yelena Mandelshtam,
\hfill  {\tt yelenam@berkeley.edu}
\vspace{-4pt}
\newpage 

\appendix

\section{Tables of (discrete) Riemann matrices}\label{appendix}
In the following tables, $n$ denotes the level of approximation as given in Definition~\ref{def:levels}. The code is always run in Matlab with run time given in seconds.
\subsection{\texorpdfstring{$L$}{L} shape tables}
\begin{table}[htbp]
    \centering
    \begin{tabular}{|c|c|c|}\hline
         $n$ & Time & Approximation\\\hline
         $0$ & 0.02 
 & $ i \begin{pmatrix} 1.75 & -1.5 \\  -1.5 & 2.00 \end{pmatrix}$ \\\hline %[0.000000000000002 + 1.750000000000000i,-0.000000000000003 - 1.500000000000000i;-0.000000000000002 - 1.500000000000000i,0.000000000000003 + 1.999999999999999i]
         $1$ & 0.05 &  $i \begin{pmatrix}1.682276986822770 & -1.364553973645541\\ -1.364553973645541 &1.729107947291081 \end{pmatrix}$ \\\hline %[0.000000000000002 + 1.682276986822770i,-0.000000000000002 - 1.364553973645541i;-0.000000000000002 - 1.364553973645541i,0.000000000000001 + 1.729107947291081i]
         $2$ &0.37& $i \begin{pmatrix} 1.670169914926280&-1.340339829852565 \\ -1.340339829852566 &1.680679659705133\end{pmatrix}$ \\\hline % [0.000000000000000 + 1.670169914926280i,-0.000000000000001 - 1.340339829852565i;0.000000000000001 - 1.340339829852566i,0.000000000000000 + 1.680679659705133i]
         $3$ &3.92 &  $ i \begin{pmatrix} 1.667472042082942 &- 1.334944084165891 \\  - 1.334944084165893 &1.669888168331791 \end{pmatrix}$ \\\hline %[0.000000000000001 + 1.667472042082942i,0.000000000000001 - 1.334944084165891i;0.000000000000001 - 1.334944084165893i,-0.000000000000002 + 1.669888168331791i]
         $4$& 28.23 & $ i \begin{pmatrix}1.666852605322711 & -1.333705210645449\\ -1.333705210645455 &1.66741042129092 \end{pmatrix}$ \\\hline %[0.000000000000001 + 1.666852605322711i,-0.000000000000004 - 1.333705210645449i;-0.000000000000003 - 1.333705210645455i,0.000000000000004 + 1.667410421290920i]
         $5$ &255.10 & $ i \begin{pmatrix}1.666709630962870 &-1.333419261925784 \\ - 1.333419261925776 & 1.666838523851582\end{pmatrix}$ \\\hline %[0.000000000000022 + 1.666709630962870i,-0.000000000000022 - 1.333419261925784i;-0.000000000000008 - 1.333419261925776i,0.000000000000011 + 1.666838523851582i]
          $6$ & 2333.12 &  $ i \begin{pmatrix} 1.666676596082551 & - 1.333353192165260\\ -1.33335319216523 &1.666706384330567 \end{pmatrix}$ \\\hline %[0.000000000000052 + 1.666676596082551i,-0.000000000000057 - 1.333353192165260i;-0.000000000000008 - 1.333353192165230i,0.000000000000025 + 1.666706384330567i]7666.3700 seconds
            $7$ & 22786.59 &  $ i \begin{pmatrix}1.666668961530435 & - 1.333337923061337\\ -1.333337923061278 &1.666675846122862 \end{pmatrix}$ \\\hline\hline %[0.000000000000158 + 1.666668961530435i,-0.000000000000170 - 1.333337923061337i;-0.000000000000011 - 1.333337923061278i,0.000000000000061 + 1.666675846122862i] A note on level 7: It took less than 12 hours, despite saying it took 26 hours =96756.3700 seconds in the cpu time.
         $\infty$ & $i \begin{pmatrix} \frac{5}{3} & -\frac{4}{3} \\ -\frac{4}{3} & \frac{5}{3} \end{pmatrix}$& $i \begin{pmatrix} 1.66666666667  & -1.33333333333  \\ -1.33333333333  &1.66666666667 \end{pmatrix}$ \\\hline 
    \end{tabular}
    \caption{Given the $L$ shape as in Figure~\ref{fig:LId} for $\lambda = 2$, the table gives the successive approximations with the bottom row representing the Riemann matrix in the limit. Since we expect the real part to be zero, we only write down the imaginary parts of the matrix level. The real parts are on the order of at worst $10^{-14}$.}
    \label{tab:newlam2}
\end{table}
\begin{table}[htbp]
    \centering
    \begin{tabular}{|c|c|c|c|}\hline
        Fraction & Tolerance & Time & $0$ level approximation \\\hline
        $\frac{15}{11}$& $1e-2$ & 0.33
 & $  i\begin{pmatrix} 1.155267361944555 & -0.582252607292078 \\  -0.582252607292077 & 1.183447277345288 \end{pmatrix}$ \\\hline %[0.000000000000001 + 1.155267361944555i,-0.000000000000002 - 0.582252607292078i;-0.000000000000001 - 0.582252607292077i,0.000000000000001 + 1.183447277345288i]
        $\frac{56}{41}$& $1e-3$ &4.04&  $ i\begin{pmatrix}1.15495696004714& -0.578505984176023\\ -0.57850598417602 &1.159755674257178 \end{pmatrix}$ \\\hline %[0.000000000000001 + 1.154956960047145i,-0.000000000000001 - 0.578505984176029i;0.000000000000001 - 0.578505984176028i,0.000000000000000 + 1.159755674257178i]
         $\frac{209}{153}$& $1e-4$ &56.96  &  $ i\begin{pmatrix} 1.154756293461396& -0.577572595239909\\- 0.577572595239901 & 1.155583435806089\end{pmatrix}$ \\\hline %[0.000000000000004 + 1.154756293461396i,-0.000000000000005 - 0.577572595239909i;0.000000000000003 - 0.577572595239901i,0.000000000000004 + 1.155583435806089i]
         $\frac{780}{571}$& $1e-5$ & 758.43&  $ i\begin{pmatrix}1.154710996247692 & -0.577390320924590\\ - 0.577390320924568 & 1.154853829287965\end{pmatrix}$ \\\hline %[0.000000000000016 + 1.154710996247692i,-0.000000000000014 - 0.577390320924590i;0.000000000000024 - 0.577390320924568i,-0.000000000000003 + 1.154853829287965i]
         $\frac{780}{571}$& $1e-6$ & 758.43 &  $ i\begin{pmatrix}1.154710996247692 &-0.577390320924590 \\-0.577390320924568 & 1.154853829287965\end{pmatrix}$ \\\hline %[0.000000000000016 + 1.154710996247692i,-0.000000000000014 - 0.577390320924590i;0.000000000000024 - 0.577390320924568i,-0.000000000000003 + 1.154853829287965i]
         $\frac{2911}{2131}$& $1e-7$ & 744.47 &  $ i\begin{pmatrix}1.15471099624769 & -0.577390320924590\\ -0.577390320924568&1.154853829287965 \end{pmatrix}$ \\\hline %[0.000000000000016 + 1.154710996247692i,-0.000000000000014 - 0.577390320924590i;0.000000000000024 - 0.577390320924568i,-0.000000000000003 + 1.154853829287965i]
         $\frac{10864}{7953}$& $1e-8$ & 774.37 &  $ i\begin{pmatrix} 1.154710996247692 &-0.57739032092459 \\-0.577390320924568 & 1.154853829287965\end{pmatrix}$ \\\hline %[0.000000000000016 + 1.154710996247692i,-0.000000000000014 - 0.577390320924590i;0.000000000000024 - 0.577390320924568i,-0.000000000000003 + 1.154853829287965i]
          $\frac{40545}{29681}$& $1e-9$ &  871.63&  $ i\begin{pmatrix} 1.154710996247692 &-0.577390320924590 \\-0.577390320924568 & 1.154853829287965\end{pmatrix}$ \\\hline\hline %[0.000000000000016 + 1.154710996247692i,-0.000000000000014 - 0.577390320924590i;0.000000000000024 - 0.577390320924568i,-0.000000000000003 + 1.154853829287965i]
        &  Exact: &$ \frac{i}{\sqrt{3}}\begin{pmatrix}2 & -1 \\ -1 & 2 \end{pmatrix}$  &  $i \begin{pmatrix}1.15470053838 & -0.57735026919\\ -0.57735026919 &1.15470053838 \end{pmatrix}$ \\\hline 
    \end{tabular}
    \caption{This table gives the successive approximations representing $\tau_\lambda$ for $\lambda = \frac{1+\sqrt{3}}{2}.$ Since we expect the real part to be zero, we only keep track of the exponential parts of the real term which are on the order of $10^{-14}$. We run a 0 level approximation, with the finer square approximations coming from increasing the tolerance according to the continued fraction expansion in Matlab.}
    \label{tab:lamsqrt}
\end{table}

\begin{table}[htbp]
    \centering
    \begin{tabular}{|c|c|c|}\hline
        Level &  Time & Approximation \\\hline
         $0$ & 774.37 &  $ i\begin{pmatrix} 1.154710996247692 &-0.57739032092459 \\-0.577390320924568 & 1.154853829287965\end{pmatrix}$ \\\hline %[0.000000000000016 + 1.154710996247692i,-0.000000000000014 - 0.577390320924590i;0.000000000000024 - 0.577390320924568i,-0.000000000000003 + 1.154853829287965i]2528.9500 seconds
          $1$ & 7290.33 &  $ i\begin{pmatrix}1.154702501426855 &- 0.577358617765855\\ - 0.577358617765802 &1.154735511279386\end{pmatrix}$ \\\hline %[0.000000000000038 + 1.154702501426855i,-0.000000000000039 - 0.577358617765855i;0.000000000000073 - 0.577358617765802i,-0.000000000000008 + 1.154735511279386i]30891.2800
          $2$ & 72388.83 &  $ i\begin{pmatrix}1.154700538230285 &- 0.577351291004541\\- 0.577351291004404 &1.154708167385750\end{pmatrix}$ \\\hline\hline %[0.000000000000151 + 1.154700538230285i,-0.000000000000121 - 0.577351291004541i;0.000000000000225 - 0.577351291004404i,-0.000000000000051 + 1.154708167385750i]438366.8600  but 24 hours real time
          $\infty$ &$\frac{10864}{7953}$ &  $i \begin{pmatrix}1.15470053534 & -0.5773502631\\  -0.5773502631&1.15470053534 \end{pmatrix}$ \\\hline 
         $\infty$ &$ \frac{1+\sqrt{3}}{2}$  &  $i \begin{pmatrix}1.15470053838 & -0.57735026919\\ -0.57735026919 &1.15470053838 \end{pmatrix}$ \\\hline 
    \end{tabular}
    \caption{This table gives an approximation of $\lambda = \frac{1+\sqrt{3}}{2}$ by the continued fraction expansion up to a tolerance of $10^{-8}$ which is $\frac{10864}{7953}$. Since we expect the real part to be zero, we only keep track of the exponential parts of the real term, which are on the order of $10^{-14}$.}
    \label{tab:lamsqrt2}
\end{table}

\newpage

\subsection{Jenkins--Strebel Tables}

\begin{table}[htbp]
    \centering
    \begin{adjustbox}{width=\textwidth}
        \begin{tabular}{|c|c|c|}\hline
         $n$ & Time & Approximation\\\hline
         $0$ & 0.01 & $\begin{pmatrix} i & 0 \\  0 & i \end{pmatrix}$ \\\hline 
         $1$ & 0.01 & {\footnotesize $\begin{pmatrix}-0.162162162162162+0.972972972972974i&-0.162162162162161-0.0270270270270267i\\
         -0.162162162162163-0.0270270270270272i&-0.162162162162162+0.972972972972973i\end{pmatrix}$}
         \\\hline 
         $2$ & 0.01& {\footnotesize$\begin{pmatrix} -0.181145110935354+0.966032669224184i&-0.181145110935353-0.0339673307758142i\\
-0.181145110935355-0.0339673307758161i&-0.181145110935356+0.966032669224183i\end{pmatrix}$}
         \\\hline
         $3$ & 0.88 & {\footnotesize $\begin{pmatrix}-0.183154151609461+0.96524676973432i&-0.183154151609459-0.0347532302656807i\\
-0.183154151609457-0.0347532302656803i&-0.183154151609459+0.965246769734321i \end{pmatrix}$}
         \\\hline 
         $4$ & 6.07 & {\footnotesize $\begin{pmatrix} -0.18337643045845+0.965159203662908i&-0.183376430458463-0.0348407963370877i\\
-0.183376430458457-0.0348407963370859i&-0.183376430458452+0.965159203662913i \end{pmatrix}$}
         \\\hline
         $5$ & 61.02 & {\footnotesize $\begin{pmatrix}-0.183401116934979+0.965149470930576i&-0.183401116935002-0.0348505290694085i\\
-0.183401116934996-0.0348505290694234i&-0.183401116934987+0.965149470930574i
 \end{pmatrix}$}
         \\\hline
         $6$ & 539.43  & {\footnotesize $\begin{pmatrix} -0.183403859739504+0.965148389476543i&-0.183403859739557-0.0348516105233965i\\
-0.183403859739527-0.0348516105234329i&-0.183403859739509+0.96514838947655i
\end{pmatrix}$}
         \\\hline
         $7$ & 3969.65 & {\footnotesize $\begin{pmatrix}-0.183404164493817+0.965148269314438i&-0.183404164494-0.0348517306853213i\\
-0.183404164493896-0.0348517306854744i&-0.183404164493843+0.965148269314462i \end{pmatrix}$}
         \\\hline
    \end{tabular}
    \end{adjustbox}
    \caption{The table gives the successive approximations of the Riemann matrix for $J_2(1,1)$, from the family of the Jenkins--Strebel differential of genus 2.}
    \label{tab:JS2}
\end{table}

\begin{table}[htbp]
    \centering
    \begin{adjustbox}{width=\textwidth}
    \begin{tabular}{|c|c|c|}\hline
         $(\lambda,\mu)$ & $n$ & Approximation\\\hline
         $(3, 1)$ &6  & {\footnotesize $\begin{pmatrix}0.391790821341348+0.799932604843609i&-0.408209178658656+0.199932604843692i\\
-0.40820917865866+0.199932604843639i&0.391790821341349+0.799932604843633i \end{pmatrix}$ }\\\hline 
         $(3/5, 1)$ & 5 & {\footnotesize $\begin{pmatrix}-0.552067431042614+0.828028726521117i&-0.0814791957486144-0.0543242146552902i\\
-0.0814791957485449-0.0543242146553346i&-0.552067431042622+0.828028726521113i\end{pmatrix}$}
         \\\hline 
         $(5/7, 1) $ & 4 & {\footnotesize$\begin{pmatrix} -0.435881709806809+0.891395736738509i&-0.111557385482526-0.0545502092073928i\\
-0.111557385482489-0.0545502092074161i&-0.435881709806802+0.891395736738513i
\end{pmatrix}$}
         \\\hline
         $(1, 3)$ & 6 & {\footnotesize $\begin{pmatrix}-0.219850572165274+0.369151081713766i&-0.0705500626168651-0.0418472644071332i\\
-0.0705500626168283-0.0418472644071658i&-0.182947389640796+0.986855533802479i \end{pmatrix}$}
         \\\hline 
         $(1, 3/5)$ & 5 & {\footnotesize $\begin{pmatrix} 0.0140501297786541+1.24511271092546i&-0.237782438362241-0.00109980530376883i\\
-0.237782438362194-0.00109980530385763i&-0.190005989628958+0.951563668084441i
\end{pmatrix}$}
         \\\hline
         $(1,5/7)$ & 4 & {\footnotesize $\begin{pmatrix}-0.0638935382171857+1.16061629684551i&-0.220568218021264-0.0143136795919605i\\
-0.220568218021245-0.0143136795919968i&-0.187013661219342+0.956221386346668i

 \end{pmatrix}$}
         \\\hline
         $(3/5, 5/7)$ & 3  & {\footnotesize $\begin{pmatrix}-0.553835613968335+1.13666373978695i&-0.110172300923-0.054295316550795i\\
-0.110172300922962-0.0542953165508437i&-0.550771147382289+0.825576938910713i

\end{pmatrix}$}
         \\\hline
    \end{tabular}
    \end{adjustbox}
    \caption{The table gives some approximations in genus 2 for $J_2(\lambda, \mu)$ where we chose the level $n$ so that the number of subdivided squares is approximately the same for each case.}
    \label{tab:JS2New}
\end{table}

\begin{table}[htbp]
    \centering
    \begin{adjustbox}{width=\textwidth}
    \begin{tabular}{|c|c|c|}\hline
         $n$ & Time & Approximation\\\hline
         $0$ & 0.01 & {\footnotesize $\begin{pmatrix} i & 0 & 0 \\  0 & i & 0 \\ 0 & 0 & 2i \end{pmatrix}$ }\\\hline  %[0.000000000000000 + 1.000000000000000i,-3.967686959186508e-16 - 2.251190959399337e-16i,1.223765959959467e-15 + 1.832448051766383e-16i;2.300213461227279e-16 + 4.202630805977024e-18i,0.000000000000000 + 1.000000000000001i,3.211959660691558e-16 - 9.266984595385614e-17i;-4.251866168899434e-16 + 7.850462293418875e-17i,-9.546942324784189e-16 - 4.960921633399546e-16i,-0.000000000000000 + 2.000000000000000i]
         $1$ & 0.01 & {\footnotesize $\begin{pmatrix}-0.163636364 + 0.972727273i &0.001474201 + 0.000245700i&-0.162162162 - 0.027027027i\\
         0.001474201 + 0.000245700i&-0.163636364 + 0.972727273i&-0.162162162 - 0.027027027i\\
         -0.162162162 - 0.027027027i&-0.162162162 - 0.027027027i&-0.324324324 + 1.945945946i\end{pmatrix}$}
         \\\hline%[-0.163636363636365 + 0.972727272727273i,0.001474201474202 + 0.000245700245701i,-0.162162162162162 - 0.027027027027028i;0.001474201474201 + 0.000245700245701i,-0.163636363636363 + 0.972727272727273i,-0.162162162162163 - 0.027027027027027i;-0.162162162162162 - 0.027027027027027i,-0.162162162162162 - 0.027027027027028i,-0.324324324324325 + 1.945945945945948i]
         $2$ & 0.01& {\footnotesize$\begin{pmatrix} -0.181890640 + 0.965429514i& 0.000745530+ 0.000603155i &-0.181145111 - 0.033967331i\\
         0.000745530 +0.000603155i &-0.181890640 + 0.965429514i&-0.181145111 - 0.033967331i
         \\-0.181145111 - 0.033967331i &-0.181145111 - 0.033967331i&-0.362290222 + 1.932065338i\end{pmatrix}$}
         \\\hline%[-0.181890640468725 + 0.965429513968413i,7.455295333670207e-04 + 6.031552557691040e-04i,-0.181145110935355 - 0.033967330775817i;7.455295333701896e-04 + 6.031552557696340e-04i,-0.181890640468726 + 0.965429513968414i,-0.181145110935355 - 0.033967330775816i;-0.181145110935356 - 0.033967330775816i,-0.181145110935352 - 0.033967330775816i,-0.362290221870711 + 1.932065338448367i]
         $3$ & 0.27 & {\footnotesize $\begin{pmatrix} -0.183837862 + 0.964626792i &0.000683710 + 0.000619978i& -0.183154152 - 0.034753230i\\
         0.006837101 + 0.000619978i&-0.183837862 + 0.964626792i&-0.183154152 - 0.034753230i\\-0.183154152 - 0.034753230i&-0.183154152 - 0.034753230i&-0.366308303 + 1.930493539i \end{pmatrix}$}
         \\\hline %[-0.183837861697396 + 0.964626792088131i,6.837100879370881e-04 + 6.199776461841601e-04i,-0.183154151609459 - 0.034753230265678i;6.837100879385080e-04 + 6.199776461851007e-04i,-0.183837861697399 + 0.964626792088131i,-0.183154151609457 - 0.034753230265678i;-0.183154151609458 - 0.034753230265676i,-0.183154151609455 - 0.034753230265678i,-0.366308303218918 + 1.930493539468640i]
         $4$ & 2.65 & {\footnotesize $\begin{pmatrix} -0.184053419 + 0.964537638i&0.000676988 + 0.000621565i&-0.183376430 - 0.034840796i\\
         0.000676988 + 0.000621565i&-0.1840534189 + 0.964537638i&-0.183376430 - 0.034840796i\\-0.183376430 - 0.034840796i&-0.183376430 - 0.034840796i&-0.366752861 + 1.930318407i\end{pmatrix}$}
         \\\hline%[-0.184053418668156 + 0.964537638481375i,6.769882096884989e-04 + 6.215651815268743e-04i,-0.183376430458449 - 0.034840796337070i;6.769882096871773e-04 + 6.215651815311827e-04i,-0.184053418668149 + 0.964537638481380i,-0.183376430458450 - 0.034840796337084i;-0.183376430458454 - 0.034840796337087i,-0.183376430458453 - 0.034840796337084i,-0.366752860916925 + 1.930318407325820i]
         $5$ &38.37 & {\footnotesize $\begin{pmatrix} -0.184077360 + 0.964527733i&0.000676243 + 0.000621738i&-0.183401117 - 0.034850529i\\
         0.000676243 + 0.000621738i&-0.184077360 + 0.964527733i&-0.183401117 - 0.034850529i\\-0.183401117 - 0.034850529i&-0.183401117 - 0.034850529i&-0.366802234 + 1.930298942i]\end{pmatrix}$}
         \\\hline%[-0.184077359904808 + 0.964527732838907i,6.762429698084880e-04 + 6.217380916470542e-04i,-0.183401116934984 - 0.034850529069379i;6.762429698119069e-04 + 6.217380916672136e-04i,-0.184077359904814 + 0.964527732838905i,-0.183401116934993 - 0.034850529069402i;-0.183401116934984 - 0.034850529069423i,-0.183401116934981 - 0.034850529069402i,-0.366802233869983 + 1.930298941861125i]
         $6$ & 398.04 & {\footnotesize $\begin{pmatrix}-0.1840800120 + 0.964526632i&0.000676160 + 0.000621757i&-0.183403860 - 0.034851611i\\
         0.000676160 + 0.000621757i&-0.1840800120 + 0.964526632i&-0.183403860 - 0.034851611i\\
         -0.183403860 - 0.034851611 i&-0.183403860 - 0.034851611i&-0.366807719 + 1.930296779i  \end{pmatrix}$}
         \\\hline%[-0.184080019924894 + 0.964526632215650i,6.761601853610258e-04 + 6.217572608186087e-04i,-0.183403859739527 - 0.034851610523303i;6.761601853490478e-04 + 6.217572608741883e-04i,-0.184080019924889 + 0.964526632215660i,-0.183403859739512 - 0.034851610523388i;-0.183403859739512 - 0.034851610523425i,-0.183403859739518 - 0.034851610523377i,-0.366807719479065 + 1.930296778953043i]
         $7$ & 3429.78 & {\footnotesize $\begin{pmatrix}-0.184080315 + 0.96452651i & 0.000676151 + 0.000621759i& -0.183404164 - 0.034851731i\\
         0.000676151 + 0.000621759i& -0.184080315 + 0.96452651i& -0.183404164 - 0.034851731i\\
         -0.183404164 - 0.034851731i& -0.183404164 - 0.034851731i & -0.366808329 + 1.930296539i  \end{pmatrix}$}
         \\\hline %[-0.184080315481226 + 0.964526509924181i,6.761509873194876e-04 + 6.217593900954986e-04i,-0.183404164493888 - 0.034851730685150i;6.761509873245022e-04 + 6.217593902530771e-04i,-0.184080315481242 + 0.964526509924194i,-0.183404164493867 - 0.034851730685331i;-0.183404164493892 - 0.034851730685467i,-0.183404164493862 - 0.034851730685328i,-0.366808328987806 + 1.930296538628769i]
    \end{tabular}
    \end{adjustbox}
    \caption{The table gives the successive approximations of the Riemann matrix for $J_3(1,1)$, a specific surface from the family of the Jenkins--Strebel differential of genus 3. Results are rounded to 9 decimal places.}
    \label{tab:JS3}
\end{table}

\begin{table}[htbp]
    \centering
     \begin{adjustbox}{width=\textwidth}
    \begin{tabular}{|c|c|c|}\hline
         $n$ & Time & Approximation\\\hline
         $0$ & 0.01 & {\footnotesize $\begin{pmatrix} i & 0 & 0  &0\\  0 & i & 0 &0 \\ 0 & 0 & i &0\\ 0 & 0 & 0 & 3i \end{pmatrix}$ }\\\hline  %See saved .txt files for exact values
         $1$ & 0.01 & {\footnotesize $\begin{pmatrix}-0.163639+0.972727i & 0.000739+0.000123i & 0.00073859+0.000123i&-0.162162-0.027027i \\
0.000739+0.000123i&-0.163639+0.972727i & 0.000739+0.000123i&-0.162162-0.027027i \\
0.000739+0.000123i&0.000739+0.000123i& -0.163639+0.972727i & -0.162162-0.027027i\\
-0.162162-0.027027i & -0.162162-0.027027i& -0.162162-0.027027i& -0.486486 +2.918919i
\end{pmatrix}$}
         \\\hline%
         $2$ & 0.02& {\footnotesize$\begin{pmatrix} -0.181893+0.96543i&0.000374+0.000301i& 0.000374+0.000301i& -0.181145-0.033967i \\
0.000374+0.000301i&-0.181893+0.965430i&0.000374+0.000301i& -0.181145-0.033967i \\
0.000374+0.000301i &0.000374+0.000301i &-0.181893+0.96543i& -0.181145-0.033967i \\
-0.181145-0.033967i&-0.181145-0.033967i&-0.181145-0.033967i&-0.543435+2.898098i
\end{pmatrix}$}
         \\\hline%
         $3$ & .48 & {\footnotesize $\begin{pmatrix}-0.18384+0.964628i& 0.000343+0.00031i& 0.000343+0.00031i&-0.183154-0.034753i \\
0.000343+0.00031i& -0.18384+0.964628i&0.000343+0.00031i& -0.183154-0.034753i \\
0.000343+0.00031i& 0.000343+0.00031i& -0.18384+0.964628i &-0.183154-0.034753i\\
-0.183154-0.034753i&-0.183154-0.034753i& -0.183154-0.034753i& -0.549462+2.89574i
 \end{pmatrix}$}
         \\\hline 
         $4$ & 5.70 & {\footnotesize $\begin{pmatrix}-0.18384+0.964628i& 0.000343+0.00031i& 0.000343+0.00031i&-0.183154-0.034753i\\
0.0003437+0.00031i& -0.18384+0.964628i& 0.000343+0.00031i& -0.183154-0.034753i \\
0.000343+0.00031i& 0.000343+0.00031i&-0.18384+0.964628i&-0.183154-0.034753i \\
-0.183154-0.034753i& -0.183154-0.034753i& -0.183154-0.034753i& -0.549462+2.89574i
 \end{pmatrix}$}
         \\\hline%
         $5$ & 89.00 & {\footnotesize $\begin{pmatrix} -0.184079+0.964529i &0.000339+0.000310i&0.000339+0.000310i&-0.183401-0.034851i \\
0.000339+0.000310i&-0.184079+0.964529i&0.000339+0.0003103i&-0.183401-0.034851i\\
0.000339+0.000310i&0.000339+0.000310i&-0.184079+0.964529i&-0.183401-0.034851i\\
-0.183401-0.0348518i&-0.183401-0.034851i&-0.183401-0.0348512i&-0.550203+2.895448i\end{pmatrix}$} 
         \\\hline%
         $6$ &  846.84 & {\footnotesize $\begin{pmatrix} -0.184082+0.9645278i&0.0003396+0.000310i&0.000339+0.000310i&-0.183404-0.034852i\\
0.000339+0.000310i&-0.184082+0.964527i&0.000339+0.000310i&-0.183404-0.034852i\\
0.000339+0.000310i&0.000339+0.000310i&-0.184082+0.964527i&-0.183404-0.034852i\\
-0.183404-0.034852i&-0.1834044-0.034852i&-0.183404-0.0348512i&-0.550212+2.895445i
 \end{pmatrix}$}
         \\\hline%[]
         $7$ & 7997.04 & {\footnotesize $\begin{pmatrix} -0.184082+0.964527i&0.000339+0.000310i&0.000339+0.000310i&-0.1834043-0.03485125i\\
0.000339+0.000310i&-0.184082+0.964527i&0.000339+0.0003105i&-0.1834041-0.0348512i\\
0.000339+0.000310i&0.000339+0.000310i&-0.184082+0.964527i&-0.183404-0.034852i\\
-0.183404-0.034852i&-0.183404-0.034852i&-0.183404-0.034852i&-0.550212+2.895445i \end{pmatrix}$}
         \\\hline %-0.184082280391336+0.964527352840432i,0.000339057948694669+0.000310458236828142i,0.000339057948610197+0.000310458236672205i,-0.183404164493683-0.0348517306848045i
%0.000339057948721604+0.000310458236985237i,-0.184082280391383+0.964527352840458i,0.000339057948558669+0.00031045823671365i,-0.183404164493691-0.0348517306850943i
%0.000339057948725268+0.000310458236989523i,0.000339057948727943+0.000310458236994094i,-0.184082280391416+0.964527352840407i,-0.183404164493851-0.03485173068532i
%-0.183404164493886-0.0348517306854523i,-0.183404164493833-0.0348517306853228i,-0.183404164493574-0.0348517306849086i,-0.550212493481976+2.89544480794274i
    \end{tabular}
    \end{adjustbox}
    \caption{The table gives the successive approximations of the Riemann matrix for $J_4(1,1)$ from the family of the Jenkins--Strebel differential of genus 4. Results rounded to 6 places.}
    \label{tab:JS4}
\end{table}
\begin{table}[htbp]
    \centering
    \footnotesize
    \begin{adjustbox}{width=\textwidth}
    \begin{tabular}{|c|c|c|}\hline
         $n$ & Time & Approximation\\\hline
         $0$ & 0.01& {\footnotesize $\begin{pmatrix} i & 0 & 0& 0 & 0\\
         0&i&0&0&0\\
         0&0&i&0&0\\
         0&0&0&i&0\\
         0&0&0&0&4i\end{pmatrix}$ }\\\hline  %
         $1$ &0.01 & { \scriptsize $ \begin{pmatrix}
         -0.163639+0.972727i&0.000737+0.000123i &0.000003&0.000737+0.000123i&-0.162162-0.027027i\\
        0.000737+0.000123i & -0.163639+0.972727i & 0.000737+0.000123i & 0.000003 &-0.162162-0.027027i\\
        0.000003 & 0.000737 +0.000123i & -0.163639+0.972727i &0.000737+0.000123i&-0.162162-0.027027i \\
        0.000737+0.000123i & 0.000003 & 0.000737+0.000123i&-0.163639+0.972727i&-0.162162 -0.027027i \\
        -0.162162-0.027027i & -0.162162-0.027027i& -0.162162 -0.027027i &-0.162162 -0.027027i&-0.648649+3.891892i
\end{pmatrix}$}
         \\\hline%-0.163639356942681+0.972726773842887i;0.00073710073710113+0.000122850122851381i\;2.99330631599e-06+4.98884386563996e-07i;0.000737100737100866+0.000122850122849\674i;-0.162162162162163-0.0270270270270269i0.000737100737100412+0.000122850122850141i;-0.163639356942678+0.972726773842888\i;0.000737100737100791+0.000122850122850944i;2.99330631579199e-06+4.98884385612\296e-07i;-0.162162162162162-0.0270270270270275i2.99330631526945e-06+4.98884385776819e-07i;0.000737100737101251+0.0001228501228\50658i;-0.163639356942678+0.972726773842887i;0.000737100737100746+0.00012285012\2849867i;-0.162162162162162-0.0270270270270284i0.000737100737100547+0.000122850122850453i;2.99330631585006e-06+4.9888438619953\4e-07i;0.00073710073710134+0.000122850122850118i;-0.16363935694268+0.9727267738\42886i;-0.162162162162162-0.0270270270270282i-0.162162162162161-0.0270270270270265i;-0.162162162162163-0.0270270270270286i;-\0.162162162162164-0.0270270270270292i;-0.162162162162163-0.0270270270270262i;-0\.648648648648648+3.89189189189189i
         $2$ & 0.03& {\scriptsize$\begin{pmatrix} 
-0.181893+0.96543i & 0.000373+0.000302i & 0.000002-0.000001i & 0.0003731+0.000303i & -0.181145-0.033967i \\
0.000373+0.000302i & -0.181893+0.96543 i & 0.000373+0.000302i& 0.000002 -0.000001i&-0.181145-0.033967i\\
0.000002 -0.000001 i& 0.000373 +0.000302i& -0.181893+0.965430i& 0.000373+0.000302i& -0.181145-0.033967i \\
0.000373+0.000302i& 0.000002- 0.000001 i& 0.000373+0.000302i &-0.1818938+0.96543i &-0.181145-0.033967i\\
-0.181145-0.033967i&-0.181145-0.033967i&-0.181145-0.033967i& -0.181145-0.033967i&-0.72458+3.864131i
\end{pmatrix}$}
         \\\hline%-0.181892971274115+0.965430427986574i;0.000372764766683059+0.000301577627889025i;2.33080538872178e-06-9.14018161716124e-07i;0.00037276476668281+0.000301577627882752i;-0.181145110935352-0.0339673307758151i
%0.000372764766686368+0.000301577627885171i;-0.181892971274115+0.965430427986577i;0.00037276476668334+0.000301577627884816i;2.33080538767077e-06-9.14018161636973e-07i;-0.18114511093535-0.0339673307758151i
%2.33080539113734e-06-9.14018160011124e-07i;0.000372764766685401+0.000301577627884619i;-0.181892971274116+0.965430427986575i;0.000372764766682164+0.00030157762788305i;-0.181145110935351-0.0339673307758168i
%0.000372764766686001+0.000301577627884512i;2.33080539067659e-06-9.1401816023436e-07i;0.000372764766685598+0.000301577627885268i;-0.181892971274118+0.965430427986575i;-0.181145110935354-0.0339673307758173i
%-0.181145110935359-0.0339673307758171i;-0.181145110935356-0.0339673307758196i;-0.181145110935353-0.0339673307758163i;-0.181145110935346-0.033967330775811i;-0.72458044374144+3.86413067689674i
         $3$ & 0.69 & {\scriptsize $\begin{pmatrix}-0.18384+0.964628i & 0.000342+0.00031i&0.000002-0.0000001i& 0.000342+0.00031i&-0.183154-0.034753i\\
0.000342+0.00031i & -0.18384+0.964628i& 0.000342+0.00031i & 0.000002-0.000001i&-0.183154-0.034753i\\
0.000002-0.000001i & 0.000342+0.00031i & -0.18384+0.964628i& 0.000342+0.00031i&-0.183154-0.034753i\\
0.000342+0.00031i&0.000002-0.000001i&0.000342+0.00031i & -0.18384+0.964628i & -0.183154 -0.034753i \\
-0.183154-0.034753i&-0.183154-0.034753i&-0.183154-0.034753i&-0.183154-0.034753i&-0.732617+3.860987i
 \end{pmatrix}$}
         \\\hline %-0.183839866512639+0.964627648348311i;0.00034185504397025+0.000309988823090987i;2.00481524148311e-06-8.56260185730889e-07i;0.000341855043969575+0.000309988823091596i;-0.183154151609467-0.0347532302656823i
%0.000341855043968496+0.000309988823093802i;-0.183839866512641+0.964627648348314i;0.000341855043971094+0.000309988823092438i;2.00481524083544e-06-8.56260184928113e-07i;-0.183154151609456-0.0347532302656849i
%2.00481523946541e-06-8.56260183737691e-07i;0.00034185504396754+0.000309988823093254i;-0.183839866512635+0.964627648348314i;0.000341855043969073+0.000309988823092448i;-0.183154151609456-0.0347532302656827i
%%0.00034185504396836+0.000309988823094077i;2.00481523888365e-06-8.56260184982677e-07i;0.000341855043969253+0.000309988823093845i;-0.183839866512637+0.964627648348313i;-0.18315415160946-0.0347532302656796i
%-0.183154151609454-0.0347532302656784i;-0.18315415160945-0.0347532302656748i;-0.183154151609462-0.0347532302656753i;-0.183154151609458-0.0347532302656755i;-0.732616606437841+3.86098707893729i

         $4$ & 9.02& {\scriptsize $\begin{pmatrix}-0.184055+0.964538i&0.000338+0.000311i &0.000002-0.000001i &0.000338+0.000311i& -0.183376-0.034841i\\
0.000338+0.000311i&-0.184055+0.964538i &0.000338+0.000311i&0.000002-0.000001i&-0.183376-0.034841i\\
0.000001-0.000001i& 0.000338+0.000311i&-0.184055+0.964538i&0.000338+0.000311i&-0.183376-0.034841i\\
0.000338+0.000311i &0.000001-0.000001i& 0.000338+0.000311i&-0.184055+0.964538i&-0.183376-0.034841i \\
-0.183376-0.034841i& -0.183376-0.034841i&-0.183376-0.034841i& -0.183376-0.034841i&-0.733506+3.860637i
 \end{pmatrix}$}
         \\\hline%-0.184055388044625+0.96453848732149i;0.000338494104851011+0.000310782590763145i;1.9693764794645e-06-8.48840120169048e-07i;0.000338494104840922+0.000310782590760174i;-0.183376430458452-0.0348407963370686i
%0.000338494104845789+0.000310782590766094i;-0.184055388044627+0.964538487321491i;0.000338494104845164+0.000310782590760257i;1.96937647648006e-06-8.48840119498947e-07i;-0.183376430458441-0.0348407963370727i
%1.9693764813324e-06-8.48840110907718e-07i;0.000338494104845445+0.000310782590763928i;-0.184055388044631+0.964538487321489i;0.000338494104839139+0.000310782590760068i;-0.183376430458448-0.0348407963370791i
%0.000338494104843603+0.000310782590767721i;1.96937648277524e-06-8.48840114003558e-07i;0.000338494104845173+0.000310782590765338i;-0.184055388044632+0.964538487321492i;-0.183376430458455-0.0348407963370843i
%-0.183376430458452-0.0348407963370953i;-0.183376430458453-0.0348407963370766i;-0.183376430458447-0.0348407963370675i;-0.183376430458438-0.0348407963370672i;-0.733505721833839+3.86063681465162i

         $5$ & 94.25 & {\scriptsize $\begin{pmatrix} 
         -0.184079+0.964529i& 0.000338+0.000311i&0.000001-0.000001i& 0.000338+0.000311i& -0.183401-0.034851i\\
0.000338+0.000311i& -0.184079+0.964529i &0.000338+0.000311i &0.000002-0.000001i&-0.183401-0.034851i\\
0.000002-0.000001i &0.000338+0.000311i& -0.184079+0.964529i &0.000338+0.000311i &-0.183401-0.034851i\\
0.000338+0.000311i& 0.000002-0.000001i &0.000338+0.000311i&-0.184079+0.964529i&-0.183401-0.034851i\\
-0.183401-0.034851i& -0.183401-0.034851i &-0.183401-0.034851i&-0.183401-0.034851i &-0.733604+3.860598i
         \end{pmatrix}$} 
         \\\hline%-0.1840793253565+0.964528580843926i;0.000338121484904282+0.00031086904582027i;1.96545169289659e-06-8.48005047351182e-07i;0.000338121484905775+0.000310869045804277i;-0.183401116934981-0.0348505290693432i
%0.000338121484903663+0.000310869045827897i;-0.184079325356516+0.964528580843924i;0.000338121484895843+0.000310869045806968i;1.96545169841755e-06-8.48005049807256e-07i;-0.183401116934982-0.0348505290693627i
%1.96545170147512e-06-8.48005026546427e-07i;0.000338121484907019+0.000310869045832612i;-0.184079325356517+0.964528580843928i;0.000338121484901591+0.000310869045806114i;-0.18340111693498-0.0348505290693834i
%0.000338121484906413+0.000310869045831275i;1.96545170434077e-06-8.48005023278086e-07i;0.000338121484907567+0.000310869045832712i;-0.184079325356512+0.964528580843925i;-0.183401116934989-0.034850529069403i
%-0.183401116934984-0.0348505290694142i;-0.183401116934983-0.034850529069398i;-0.183401116934974-0.0348505290693721i;-0.183401116934977-0.0348505290693414i;-0.733604467739997+3.86059788372217i

         $6$ & 1051.54  & {\scriptsize $\begin{pmatrix} 
 -0.184082+0.964527i& 0.000338+0.000311i & 0.000002-0.000001i& 0.000338+0.000311i& -0.183404-0.034852i\\
0.000338+0.000311i& -0.184082+0.964527i& 0.000338+0.000311i& 0.000002-0.000001i& -0.183404-0.034852i\\
0.000002-0.000001i& 0.000338+0.000311i& -0.184082+0.964527i& 0.000338+0.000311i& -0.183404-0.034852i\\
0.000338+0.000311i& 0.000002-0.000001i& 0.000338+0.000311i& -0.184082+0.964527i& -0.183404-0.034852i\\
-0.183404-0.034852i& -0.183404-0.034852i& -0.183404-0.034852i& -0.183404-0.034852i& -0.733615+3.860594i\end{pmatrix}$}
         \\\hline%[]
         $7$ & 12739.99  & {\scriptsize $\begin{pmatrix} -0.184082+0.964527i& 0.000338+0.000311i& 0.000002-0.000001i& 0.000338+0.000311i& -0.183404-0.034852i\\
0.000338+0.000312i& -0.184082+0.964527i& 0.000338+0.000311i& 0.000002-0.000001i & -0.183404-0.034852i\\
0.000002 -0.000001i& 0.000338+0.000311i& -0.184082+0.964527i& 0.000338+0.000311i& -0.183404-0.034852i\\
0.000338+0.000311i& 0.000002 - -0.000001i& 0.000338+0.000311i& -0.184082+0.964527i& -0.183404-0.034852i\\
-0.183404-0.034852i& -0.183404-0.034852i& -0.183404-0.034852i& -0.183404-0.034852i& -0.733617+3.860593i \end{pmatrix}$}
         \\\hline 
    \end{tabular}
    \end{adjustbox}
    \caption{The table gives the successive approximations of the Riemann matrix for $J_5(1,1)$ from the family of the Jenkins--Strebel differential of genus 5. Results rounded to 6 places.}
    \label{tab:JS5}
\end{table}
\end{document}